\documentclass[a4paper,12pt]{report}

\newtheorem{proposition}{Proposition}[chapter]

\newtheorem{conjecture}{Conjecture}[chapter]

\newtheorem{definition}{Definition}[chapter]
\newtheorem{assumption}{Assumption}[chapter]
\newtheorem{remark}{Remark}[chapter]

\newtheorem{example}{Example}[chapter]

\usepackage{amsmath}
\usepackage{amssymb}
\usepackage{amsfonts}
\usepackage{graphicx}
\usepackage{mathrsfs}

\title{\bf Optimization on linear matrix inequalities\\for polynomial systems control}
\begin{document}

\author{Didier Henrion$^{1,2}$}

\footnotetext[1]{LAAS-CNRS, Univ. Toulouse, France. {\tt henrion@laas.fr}}
\footnotetext[2]{Fac. Elec. Engr., Czech Tech. Univ. Prague, Czech Rep. {\tt henrion@fel.cvut.cz}}

\date{Draft lecture notes of \today}

\maketitle

\centerline{\bf Abstract}

Many problems of systems control theory boil down to solving polynomial equations,
polynomial inequalities or polyomial differential equations. Recent advances in convex
optimization and real algebraic geometry can be combined to generate
approximate solutions in floating point arithmetic. 

In the first part of the course
we describe semidefinite programming (SDP) as an extension of linear programming (LP)
to the cone of positive semidefinite matrices.
We investigate the geometry of spectrahedra, convex sets defined by 
linear matrix inequalities (LMIs) or affine sections of the SDP cone.
We also introduce spectrahedral shadows, or lifted LMIs, obtained
by projecting affine sections of the SDP cones.
Then we review existing numerical algorithms for solving SDP problems.

In the second part of the course 
we describe several recent applications of SDP. First, we explain how to solve
polynomial optimization problems, where a real multivariate polynomial must
be optimized over a (possibly nonconvex) basic semialgebraic set. 
Second, we extend these techniques to ordinary differential equations (ODEs)
with polynomial dynamics, and the problem of trajectory optimization (analysis
of stability or performance of solutions of ODEs). Third, we conclude
this part with applications to optimal control (design of a trajectory optimal
w.r.t. a given functional).

For some of these decision and optimization problems,
it is hoped that the numerical solutions computed by SDP can be refined a posteriori
and certified rigorously with appropriate techniques.

\vspace{2em}

\centerline{\bf Disclaimer}

These lecture notes were written for a tutorial course given during
the conference ``Journ\'ees Nationales de Calcul Formel'' held at
Centre International
de Rencontres Math\'ematiques, Luminy, Marseille, France in May 2013.
They are aimed at giving an elementary and introductory account
to recent applications of semidefinite programming to the numerical
solution of decision problems involving polynomials in
systems and control theory. The main technical results are gathered in
a hopefully concise, notationally simple way, but for the sake
of conciseness and readability, they are not proved in the document.
The reader interested in mathematical rigorous comprehensive
technical proofs is referred to the papers and books listed
in the ``Notes and references'' section of each chapter.
Comments, feedback, suggestions for improvement of these
lectures notes are much welcome.

\vspace{2em}

\centerline{\bf Acknowledgement}

I am grateful to the organizers of the conference,
and especially to Guillaume Ch\`eze and Mohab Safey El Din, for giving
me the opportunity to prepare and present this material.
Many thanks to Mohamed Rasheed Hilmy Abdalmoaty, Mathieu Claeys,
Simone Naldi, Luis Rodrigues and Mohab Safey El Din for their
suggestions and remarks on the text.

\newpage

\tableofcontents

\chapter{Motivating examples}

In this introductory section we describe elementary problems
of systems control theory that can be formulated as decision and optimization
problems over polynomial equations and differential equations.

\section{Structured eigenvalue assignment}\label{eig}

We consider the problem of designing a pulsed power generator in an electrical network.
The engineering specification of the design is that a suitable resonance condition
is satisfied by the circuit so that the energy initially stored in a
number of stage capacitors is transferred in finite time to a single load capacitor
which can then store the total energy and deliver the pulse.

Mathematically speaking, the problem can be formulated as
the following structured matrix eigenvalue assignment problem.
Let $n \in {\mathbb N}$ and define the matrix
\[
B =
\left(\begin{array}{ccccc}
2 & -1 & 0 & \cdots & 0\\
-1 & 2 &  &  & \vdots\\
0 & & \ddots &  & 0\\
\vdots & & & 2 & -1 \\
0 & \cdots & 0 & -1 &  \frac{n + 1}{n}
\end{array}\right).
\]
Given positive rational numbers $a_k \in {\mathbb Q}$, $k=1,\ldots,n$,
consider the eigenvalue assignment problem
\[
\det(sI_n-B^{-1}\:\mathrm{diag}\:x) = s^n+b_1(x)s^{n-1}+\cdots+b_{n-1}(x)s+b_n(x) = \prod_{k=1}^n (s+a_k)
\]
where $x \in {\mathbb R}^n$ is a vector of unknowns
and $\mathrm{diag}\:x$ stands for the $n$-by-$n$
matrix with entries of $x$ along the diagonal.
In systems and control terminology, this is
a structured pole placement problem, and vector $x$
can be interpreted as a parametrization of a 
linear controller to be designed.
By identifying like powers of indeterminate $s$
in the above relation,
it can be formulated as a polynomial system
of equations $p_k(x)=0$, $k=1,\ldots,n$ where
\[
\begin{array}{rcl}
p_1(x) & = & b_1(x) - a_1 - \cdots - a_n \\
p_2(x) & = & b_2(x) - a_1 a_2 - a_1 a_3 - \cdots - a_{n-1} a_n \\
& \vdots \\
p_n(x) & = & b_n(x) - a_1 a_2 \cdots a_n.
\end{array}
\]
In the context of electrical generator design,
a physically relevant choice of eigenvalues is
\[
a_k = \frac{1}{(2k)^2-1}, \quad k=1,\ldots,n.
\]
For example, if $n=2$, we obtain the following system
\[
\begin{array}{rcl}
\frac{3}{4} x_1+x_2-\frac{2}{5} &=& 0 \\
\frac{1}{2} x_1 x_2-\frac{1}{45} &=& 0. \\
\end{array}
\]
More generally, we obtain a system with $n$ unknowns and $n$ polynomial equations 
of respective degrees $1,\ldots,n$
which has typically much less than $n!$ real solutions. Geometrically,
the feasibility set
\[
X = \{x \in {\mathbb R}^n \: :\: p_k(x)=0,\quad k=1,\ldots,n\}
\]
is a zero-dimensional real algebraic set of small cardinality. 
When $n=8$, say, we would like to find a point in $X$.

\section{Control law validation}\label{valid}

In aerospace engineering, the validation of control laws
is a critical step before industrialization. Generally it is
carried out by expensive time-simulations.
A very simple, but representative example, is the validation
of a control law for a one-degree-of-freedom
model of a launcher attitude control system in orbital
phase. The closed-loop system must follow a given
piecewise linear angular velocity profile. It is modeled
as a double integrator
\[
I \ddot{\theta}(t) = u(t)
\]
where $I$ is a given constant inertia, $\theta(t)$ is the angle
and $u(t)$ is the torque control. We denote
\[
x(t) = \left[\begin{array}{c} \theta(t) \\ \dot{\theta}(t) \end{array}\right]
\]
and we assume that both angle $x_1(t)$ and angular velocity $x_2(t)$
are measured, and that the torque control is given by
\[
u(x(t)) = \mathrm{sat}(F'\mathrm{dz}(x_r(t)-x(t)))
\]
where $x_r(t)$ is the given reference signal, $F \in {\mathbb R}^2$
is a given state feedback, the prime denotes transposition,
$\mathrm{sat}$ is a saturation function such that $\mathrm{sat}(y) = y$ if
$|y|\leq L$ and $\mathrm{sat}(y) = L\:\mathrm{sign}(y)$ otherwise, $\mathrm{dz}$ is
a dead-zone function such that $\mathrm{dz}(x) = 0$ if $|x_i| \leq D_i$ for some $i=1,2$
and $\mathrm{dz}(x) = 1$ otherwise. Thresholds $L>0$, $D_1>0$ and $D_2>0$ are given.

We would like to verify whether the system state $x(t)$
reaches a given subset $X_T = \{x \in {\mathbb R}^2 \: :\: x^Tx \leq \varepsilon\}$
of the deadzone region after a fixed time $T$, and for all possible initial conditions $x(0)$ chosen
in a given subset $X_0$ of the state-space, and for zero reference signals.

\section{Bolza's optimal control problem}\label{bolza}

\begin{figure}[th]
\begin{center}
\includegraphics[width=0.9\textwidth]{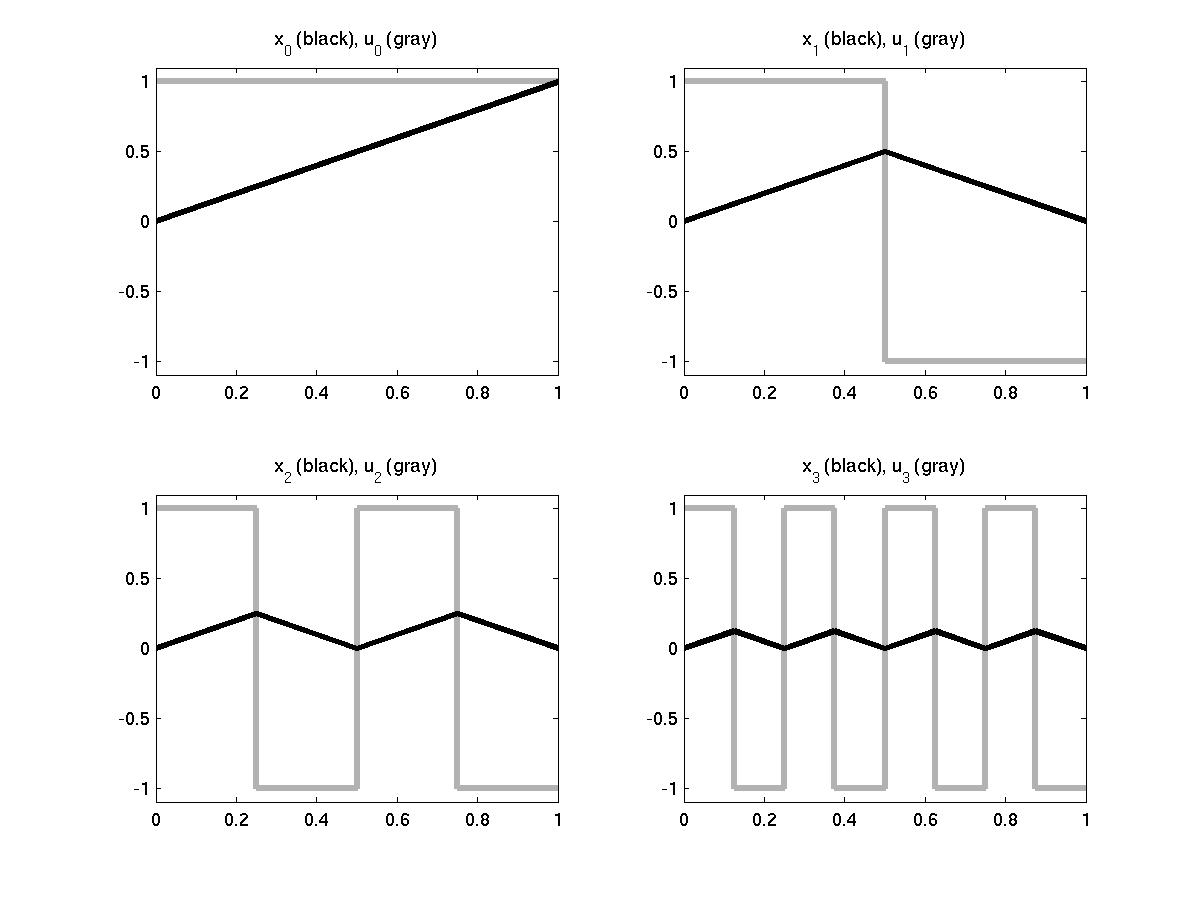}
\caption{Sequences of state trajectories and control inputs for Bolza's example.\label{fig:bolza}}
\end{center}
\end{figure}

Our last example is a classical academic problem of calculus of variations
illustrating that an optimal control problem with smooth data (infinitely differentiable
Lagrangian and dynamics, no state or input constraints) can have
a highly oscillatory optimal solution.

Consider the optimal control problem 
\[
\begin{array}{lll}
p^* = & \inf & \int_0^1 \left(x^4(t)+(u^2(t)-1)^2\right) dt\\ 
& \mathrm{s.t.} &  \dot{x}(t) = u(t), \quad t \in [0,1], \\ 
& & x(0) = 0, \quad x(1) = 0
\end{array}
\]
where the infimum is w.r.t. a Lebesgue integrable real-valued
control $u \in {\mathscr L}^1([0,1]; {\mathbb R})$.
Intuitively, the state trajectory
$x(t)$ should remain zero, and the velocity $\dot{x}=u$ should
be equal to $+1$ or $-1$, so that the nonnegative Lagrangian
$l(t,x(t),u(t)):=x^4(t)+(u^2(t)-1)^2$ remains zero, and hence the
objective function is zero, the best we can hope.
We can build a sequence of bang-bang controls $u_k(t)$
such that for each $k=0,1,2,\ldots$ the corresponding state trajectory $x_k(t)$
has a sawtooth shape, see Figure \ref{fig:bolza}. With such a sequence the
objective function tends to $\lim_{k\to\infty} \int_0^1 l(t,x_k(t),u_k(t)) dt =
\int_0^1 x^4_k(t)dt = 0$ and hence $p^*=0$. 
This infimum is however not attained with a control law $u(t)$
belonging to the space of Lebesgue integrable functions.

We would like to develop a numerical method that can deal
with such oscillation phenomena and would allow us to
construct explicitly an optimal control law.

\section{Course outline}

The objective of this document is to describe a systematic
approach to the numerical solution of these nonlinear
nonconvex decision problems. Our strategy will be
as follows:
\begin{enumerate}
\item the problem is relaxed and linearized to an LP on measures,
interpreted as the dual to an LP on continuous functions;
\item since the decision problems have polynomial data, the measure LP
is formulated as a moment LP;
\item a hierarchy of finite-dimensional LMI
relaxations is used to solve the moment LP numerically,
with guarantees of asymptotic, and sometimes finite
convergence.
\end{enumerate}

Since we do not assume that the reader is familiar
with SDP and the geometry of LMIs, the document
starts with an introductory Chapter \ref{chap:finlp}
on finite-dimensional conic programming. In Chapter
\ref{chap:polyopt}, our approach
is applied to nonconvex finite-dimensional polynomial
optimization. Finally, we conclude with Chapter
\ref{chap:ido} on nonconvex infinite-dimensional optimization
on solutions of polynomial differential equations,
and a last Chapter \ref{chap:poc} on extensions to polynomial
optimal control.


\chapter{Conic optimization}\label{chap:finlp}

\section{Convex cones}

In this section we describe linear programming over
convex cones in finite dimensional Euclidean spaces.

\begin{definition}[Convex set]
A set $K$ is convex if $x_1,x_2 \in K$ implies
$\lambda x_1 + (1-\lambda) x_2 \in K$ for all $\lambda \in [0,1]$.
\end{definition}

Geometrically speaking, a set is convex whenever the line segment
linking two of its points belongs to the set.

\begin{definition}[Convex hull]
Given a set $K$, its convex hull, denoted by $\mathrm{conv}\:K$,
is the smallest closed convex set which contains $K$.
\end{definition}

The convex hull can be expressed as
\[
\mathrm{conv}\:K:=\left\{\sum_{k=1}^N \lambda_k x_k \: : \: N \in {\mathbb N}, 
\: x_k \in K, \: \lambda_k \geq 0, \: \sum_{k=1}^N \lambda_k = 1\right\}.
\]
The convex hull of finitely many points $\mathrm{conv}\{x_1,\ldots,x_N\}$ 
is a polytope with vertices at these points.
A theorem by Carath\'eodory states that given a set $K \subset {\mathbb R}^n$,
every point of $\mathrm{conv}\:K$ can be expressed as $\sum_{k=1}^{n+1} \lambda_k x_k$
for some choice of $x_k \in K$, $\lambda_k \geq 0$, $\sum_{k=1}^{n+1} \lambda_k=1$.

\begin{definition}[Cone]
A set $K$ is a cone if $\lambda \geq 0, x \in K$ implies $\lambda x \in K$.
\end{definition}

It follows that a convex cone is a set which is invariant under
addition and multiplication by non-negative scalars.

Let us denote the scalar product of two vectors $x,y$ of ${\mathbb R}^n$ as follows:
\[
\langle x,y \rangle := x'y = \sum_{k=1}^n x_k y_k
\]
where the prime, applied to a vector or a matrix, denotes the transpose.
More generally, we use the prime to denote the dual vector space:
\begin{definition}[Dual space]\label{dualspace}
The dual of a vector space $V$ is the space $V'$ of all linear functionals on $V$.
\end{definition}

When applied to a cone, the prime denotes the dual cone:
\begin{definition}[Dual cone]\label{dualcone}
The dual of a cone $K$ is the cone
\[
K':=\{y \in {\mathbb R}^n \: :\: \langle x,y \rangle \geq 0, \:\:\forall x \in K\}.
\]
\end{definition}

Geometrically, the dual cone of $K$ is the set of all nonnegative linear
functions on $K$. 
Notice that the dual $K'$ is always a closed convex cone
and that $K''$ is the closure of the conic hull, i.e.
the smallest convex cone that contains $K$.
In particular, if $K$ is a closed convex cone, then $K''=K$.
A cone $K$ such that $K'=K$ is called self-dual.

A cone $K$ is pointed if $K \cap (-K) = \{0\}$ and solid if the interior of $K$ is not empty.
A cone which is convex, closed, pointed and solid is called a proper cone.
The dual cone of a proper cone is also a proper cone.
A proper cone $K$ induces a partial order (a binary relation that is reflexive,
antisymmetric and transitive) on the vector space: $x_1 \geq x_2$ if and only if
$x_1-x_2 \in K$.

\begin{definition}[Linear cone]
The linear cone, or positive orthant, is the set
\[
\{x \in {\mathbb R}^n \: :\: x_k \geq 0, \: k=1,\ldots,n\}.
\]
\end{definition}

\begin{definition}[Quadratic cone]
The quadratic cone, or Lorentz cone, is the set
\[
\{x \in {\mathbb R}^n \: :\: x_1 \geq \sqrt{x^2_2+\cdots+x^2_n}\}.
\]
\end{definition}

Let ${\mathbb S}^n$ denote the Euclidean space of $n$-by-$n$ symmetric matrices of
${\mathbb R}^{n\times n}$, equipped with the inner product
\[
\langle X,Y \rangle := \mathrm{trace}\:X'Y = \sum_{i=1}^n \sum_{j=1}^n x_{ij}y_{ij}
\]
defined for two matrices $X$, $Y$ with respective entries $x_{ij}$, $y_{ij}$,
$i,j=1,\ldots,n$.

\begin{definition}[Gram matrix]
Given a real quadratic form $f : {\mathbb R}^n \to {\mathbb R}$,
the (unique) matrix $X \in {\mathbb S}^n$ such that $f(y)=y'Xy$
is called the Gram matrix of $f$.
\end{definition}

\begin{definition}[Positive semidefinite matrix]
A matrix is positive semidefinite when it is the Gram matrix
of a nonnegative quadratic form.
\end{definition}

In other words, a matrix $X \in {\mathbb S}^n$ is positive semidefinite,
denoted by $X \geq 0$, if and only if  $y'Xy \geq 0, \:\forall y \in {\mathbb R}^n$
or equivalently, if and only if the minimum eigenvalue of $X$ is non-negative.
This last statement makes sense since symmetric matrices have only real eigenvalues. 

\begin{definition}[Semidefinite cone]
The semidefinite cone, or cone of positive semidefinite matrices, is the set
\[
\{X \in {\mathbb S}^n \: :\: X \geq 0\}.
\]
\end{definition}

\begin{proposition}[Self-dual cones]\label{selfdual}
The linear, quadratic and semidefinite cones are self-dual convex cones.
\end{proposition}

Note finally that if $K={\mathbb R}^n$ is interpreted as a cone, then
its dual $K'=\{0\}$ is the zero cone, which contains only
the zero vector of ${\mathbb R}^n$.

\section{Primal and dual conic problems}

{\bf Conic programming} is linear programming in a convex cone $K$: we want to
minimize a linear function over the intersection of $K$ with an affine subspace:
\begin{equation}\label{plp}
\begin{array}{rcll}
p^* & = & \inf & c'x \\
&& \mathrm{s.t.} & Ax = b \\
&&& x \in K
\end{array}
\end{equation}
where the infimum is w.r.t. a vector $x \in {\mathbb R}^n$ to be found,
and the given problem data consist of a matrix $A \in {\mathbb R}^{m\times n}$,
a vector $b \in {\mathbb R}^m$ and a vector $c \in {\mathbb R}^n$.
Note that the feasibility set $\{x \in {\mathbb R}^n \: :\: Ax = b, \: x \in K\}$
is not necessarily closed, so that in general we speak of an infimum, not
of a minimum.

If $K={\mathbb R}^n$, the whole Euclidean space, or free cone,
then problem (\ref{plp}) amounts to solving a linear system of equations.
If $K$ is the linear cone, then solving problem (\ref{plp}) is called linear programming (LP).
If $K$ is the quadratic cone, then this is called (convex) quadratic programming (QP).
If $K$ is the semidefinite cone, then this is called (linear) {\bf semidefinite programming} (SDP).

In standard mathematical programming terminology, problem (\ref{plp}) is called the primal
problem, and $p^*$ denotes its infimum. The primal conic problem has a dual conic problem:
\begin{equation}\label{dlp}
\begin{array}{rcll}
d^* & = & \sup & b'y \\
&& \mathrm{s.t.} & z=c-A'y \\
&&& z \in K'.
\end{array}
\end{equation}
Note that from Proposition \ref{selfdual}, if $K$ is the direct product
of linear, quadratic and semidefinite cones, then $K'=K$.
If $K$ contains a free cone, then the corresponding components in $K'$
are zero:  we can enforce equality constraints on
some entries in vector $z$ in dual problem (\ref{dlp}), and they
correspond to unrestricted entries in vector $x$ in primal problem (\ref{plp}).

\begin{example}
If $K$ is the direct product of a 2-dimensional free cone with a 2-dimensional
LP cone and a 2-dimensional SDP cone, then in primal problem (\ref{plp})
the constraint $x \in K \subset {\mathbb R}^7$ can be expressed entrywise as:
\[
\begin{array}{l}
x_1 \:{\rm free}\:, x_2 \:{\rm free},\\
x_3 \geq 0, \: x_4 \geq 0, \\
\left(\begin{array}{cc}
x_5 & x_6 \\ x_6 & x_7
\end{array}\right) \geq 0
\end{array}
\]
and in dual problem (\ref{dlp}) the constraint $z \in K' \subset {\mathbb R}^7$
can be expressed entrywise as:
\[
\begin{array}{l}
z_1 = 0, \: z_2 = 0, \\
z_3 \geq 0, \: z_4 \geq 0, \\
\left(\begin{array}{cc}
z_5 & z_6 \\ z_6 & z_7
\end{array}\right) \geq 0.
\end{array}
\]
\end{example}

If $K$ consists of only one semidefinite cone, primal problem (\ref{plp})
can be written as follows:
\begin{equation}\label{psdp}
\begin{array}{rcll}
p^* & = & \inf & \langle C,X \rangle \\
&& \mathrm{s.t.} & {\mathcal A}X = b \\
&&& X \geq 0
\end{array}
\end{equation}
where the given problem data consist now of a linear operator
${\mathcal A} : {\mathbb S}^n \to {\mathbb R}^m$,
a vector $b \in {\mathbb R}^m$ and a matrix $C \in {\mathbb S}^n$.
The action of operator $\mathcal A$ is described entrywise as $\langle A_k,X \rangle = b_k$,
for given matrices $A_k \in {\mathbb S}^n$, $k=1,\ldots,m$.
The adjoint or dual operator ${\mathcal A}' : ({\mathbb R}^m)'={\mathbb R}^m \to
({\mathbb S}^n)'={\mathbb S}^n$
is the unique linear map such that
$\langle {\mathcal A}'y,X \rangle = \langle y,{\mathcal A}X\rangle$
for all $X \in {\mathbb S}^n$ and $y \in {\mathbb R}^m$.
More concretely, ${\mathcal A}'y = \sum_{k=1}^m A_k y_k$.

Primal SDP problem (\ref{psdp}) has a dual SDP problem:
\begin{equation}\label{dsdp}
\begin{array}{rcll}
d^* & = & \sup & \langle b,y \rangle \\
&& \mathrm{s.t.} & C-{\mathcal A}'y \geq 0
\end{array}
\end{equation}
where the supremum is w.r.t. a vector $y \in {\mathbb R}^m$.

\begin{example}
The primal SDP problem
\[
\begin{array}{rcll}
p^* & = & \inf & x_{11}+x_{22}+x_{33} \\
&& \mathrm{s.t.} & -2x_{21} = 1 \\
&&& -2x_{31} = 1 \\
&&& -2x_{32} = 1 \\
&&& \left(\begin{array}{ccc}x_{11}&x_{21}&x_{31}\\
x_{21}&x_{22}&x_{32}\\x_{31}&x_{32}&x_{33}\end{array}\right) \geq 0
\end{array}
\]
has a dual SDP problem
\[
\begin{array}{rcll}
d^* & = & \sup & y_1+y_2+y_3 \\
&& \mathrm{s.t.} & \left(\begin{array}{ccc}
1&y_1&y_2 \\ y_1&1&y_3 \\ y_2&y_3&1\end{array}\right) \geq 0.
\end{array}
\]
Both problems share the data
\[
A_1 = -\left(\begin{array}{ccc}0&1&0\\1&0&0\\0&0&0\end{array}\right), \:
A_2 = -\left(\begin{array}{ccc}0&0&1\\0&0&0\\1&0&0\end{array}\right), \:
A_3 = -\left(\begin{array}{ccc}0&0&0\\0&0&1\\0&1&0\end{array}\right)
\]
and
\[
b = \left(\begin{array}{c}1\\1\\1\end{array}\right), \:
C = \left(\begin{array}{ccc}1&0&0\\0&1&0\\0&0&1\end{array}\right)
\]
on the 3-dimensional SDP cone $K=K'$.
\end{example}

\section{Spectrahedra and LMIs}\label{lmi}

The convex feasibility sets of problems (\ref{plp}) and (\ref{dlp})
are intersections of a convex cone with an affine subspace.
We would like to understand the geometry of these sets.
In particular, we would like to know whether a given convex set
can be modeled like this.

The most general case relevant for our purposes
is when $K$ is the direct product of semidefinite cones.
Indeed, note first that every linear cone is the
direct product of one-dimensional quadratic cones,
or equivalently, of one-dimensional semidefinite cones.
Second, note that a quadratic cone is a particular
affine section of the semidefinite cone:
\[
\{x \in {\mathbb R}^n \: :\: x_1 \geq \sqrt{x^2_2+\cdots+x^2_n}\} \:=\: 
\{x \in {\mathbb R}^n \: :\:
\left(\begin{array}{cccc}
x_1 & x_2 & \cdots & x_n \\
x_2 & x_1 & & 0 \\
\vdots & & \ddots & \vdots \\
x_n & 0 & \cdots & x_1
\end{array}\right) \geq 0\}.
\]
It follows that every set that can be represented as
an affine section of direct products of the linear and
quadratic cone can be represented as an affine section
of direct products of the semidefinite cone.
Finally, note that a direct product of semidefinite cones
can be expressed as an affine section of a single semidefinite cone, e.g.
\[
\{x \in {\mathbb R}^4 \: :\: x_1 \geq 0, \: \left(\begin{array}{cc}
x_2 & x_3 \\ x_3 & x_4\end{array}\right) \geq 0\} \:=\:
\{x \in {\mathbb R}^4 \: :\: \left(\begin{array}{ccc}
x_1 & 0 & 0 \\ 0 & x_2 & x_3 \\ 0 & x_3 & x_4
\end{array}\right) \geq 0\}.
\]
For this reason, in most of the remainder of this document,
we consider a single semidefinite cone constraint.

\begin{definition}[LMI]
A linear matrix inequality (LMI) is a constraint
\[
F_0 + \sum_{k=1}^n x_k F_k \geq 0
\]
on a vector $x \in {\mathbb R}^n$, where matrices
$F_k \in {\mathbb S}^m$, $k=0,1,\ldots,n$ are given.
\end{definition}

Note that an LMI constraint is generally nonlinear,
but it is always convex. To prove convexity, rewrite
the LMI constraint as
\[
y'\left(F_0 + \sum_{k=1}^n x_k F_k\right)y = (y'F_0y) + \sum_{k=1}^n (y'F_ky) x_k \geq 0
\]
which models infinitely many linear constraints on $x \in {\mathbb R}^n$,
parametrized by $y \in {\mathbb R}^m$.

\begin{definition}[Spectrahedron, or LMI set]
A spectrahedron is a set described by an LMI:
\[
\{x \in {\mathbb R}^n \: :\: F_0 + \sum_{k=1}^n x_k F_k \geq 0\}
\]
where matrices $F_k \in {\mathbb S}^m$, $k=0,1,\ldots,n$ are given.
\end{definition}

In other words, spectrahedra are affine sections of
the semidefinite cone, or equivalently, LMI sets.
Note that in the case where matrices $F_k$, $k=0,1,\ldots,n$ all commute (e.g.
if they are all diagonal), the LMI reduces to $m$ affine inequalities,
and the spectrahedron reduces to a polyhedron.

\begin{figure}
\begin{center}
\includegraphics[width=0.5\textwidth]{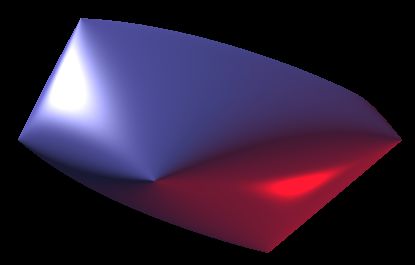}
\caption{A spectrahedron.\label{detquintic}}
\end{center}
\end{figure}

On Figure \ref{detquintic} we represent a spectrahedron in the case $n=3$
and $m=5$. We observe that its boundary is almost everywhere smooth
and curved outwards (by convexity), but it also includes vertices and
edges.

Let ${\mathbb R}[x]$ denote the ring of polynomials of the
indeterminate $x \in {\mathbb R}^n$ with real coefficients.
Given a polynomial $f \in {\mathbb R}[x]$,
we define its set of zeros, or level set, as
$\{x \in {\mathbb R}^n \: :\: f(x) = 0\}$.
We define its open superlevel set as
$\{x \in {\mathbb R}^n \: :\: f(x) > 0\}$,
and its closed superlevel set as
$\{x \in {\mathbb R}^n \: :\: f(x) \geq 0\}$.
Note that these sets are defined in ${\mathbb R}^n$,
not in ${\mathbb C}^n$, since in this document
we are mainly concerned with optimization.

\begin{definition}[Algebraic set]
An algebraic set is an intersection
of finitely many polynomial level sets.
\end{definition}

\begin{definition}[Semialgebraic set]
A semialgebraic set is a union of
finitely many intersections of finitely many
open polynomial superlevel sets.
\end{definition}

\begin{definition}[Closed basic semialgebraic set]
A closed basic semialgebraic set is an
intersection of finitely many closed polynomial superlevel sets.
\end{definition}

Now, let us denote by
\[
F(x) := F_0 + \sum_{k=1}^n x_k F_k 
\]
the affine symmetric matrix describing a spectrahedron, and build its
characteristic polynomial
\[
t \mapsto \det\:(tI_m+F(x)) = \sum_{k=0}^m f_{m-k}(x)t^k
\]
which is monic, i.e. $f_0(x)=1$. Coefficients $f_k \in {\mathbb R}[x]$,
$k=1,\ldots,m$
are multivariate polynomials called the defining polynomials
of the spectrahedron. They are elementary symmetric functions
of the eigenvalues of $F(x)$.

\begin{proposition}[Spectrahedra are closed basic semialgebraic sets]\label{csa}
A spectrahedron can be expressed as follows:
\[
\{x \in {\mathbb R}^n \: :\: F_0 + \sum_{k=1}^n x_k F_k \geq 0\} \:=\: 
\{x \in {\mathbb R}^n \: :\: f_k(x) \geq 0, \: k=1,\ldots,m\}.
\]
\end{proposition}

\begin{figure}
\begin{center}
\includegraphics[width=0.5\textwidth]{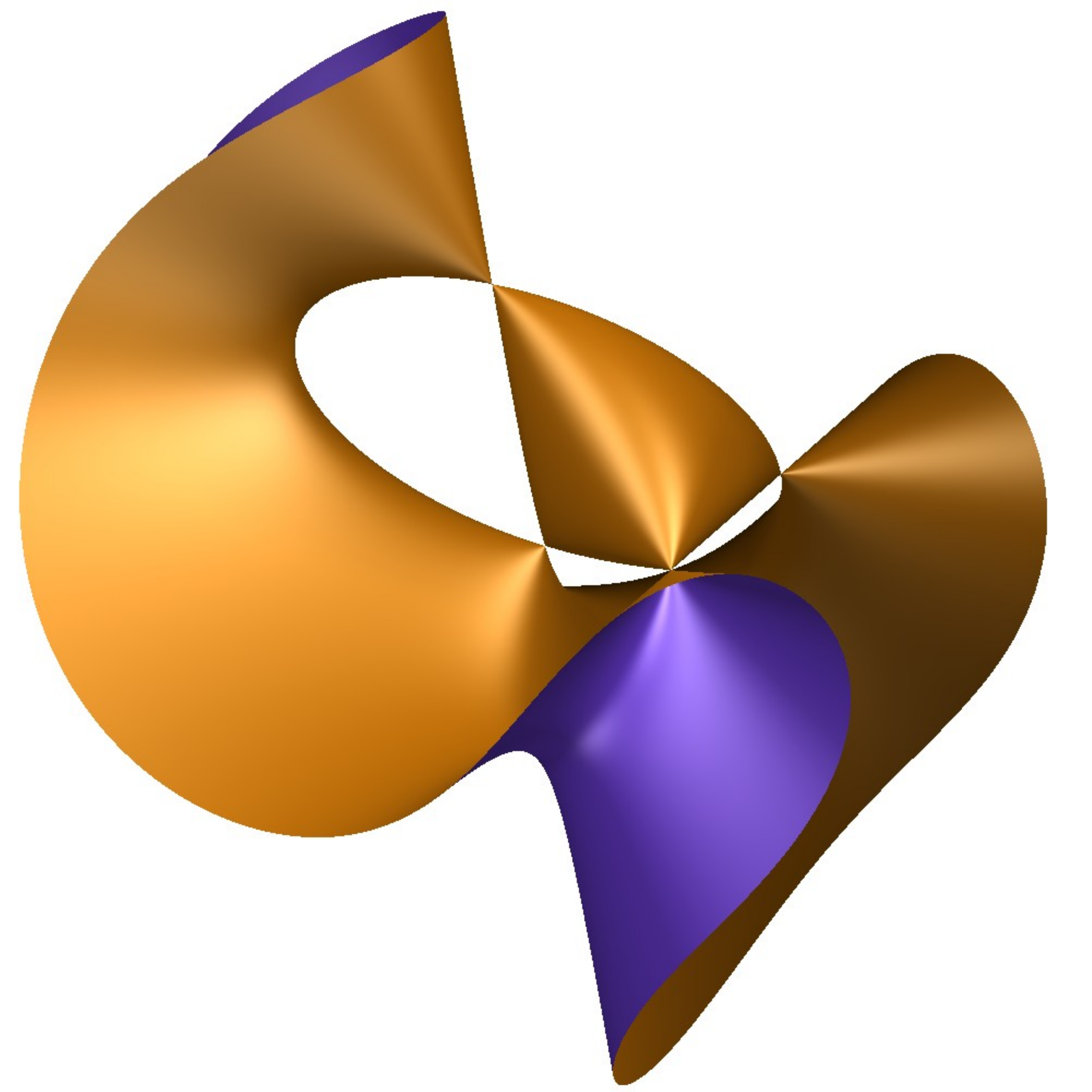}
\caption{The Cayley cubic surface and its spectrahedron.\label{cayley_cubic}}
\end{center}
\end{figure}

\begin{example}[The pillow]\label{ex:pillow}
As an elementary example, consider the pillow spectrahedron
\[
X:=\{x \in {\mathbb R}^3 \: :\: F(x) := \left(\begin{array}{ccc}
1 & x_1 & x_2 \\ x_1 & 1 & x_3 \\ x_2 & x_3 & 1
\end{array}\right) \geq 0 \}
\]
and its defining polynomials
\[
\begin{array}{rcl}
f_1(x) & = & \mathrm{trace}\:F(x) = 3, \\
f_2(x) & = & 3-x^2_1-x^2_2-x^2_3, \\
f_3(x) & = & \det F(x) = 1+2x_1x_2x_3-x^2_1-x^2_2-x^2_3,
\end{array}
\]
On Figure \ref{cayley_cubic} we represent the Cayley cubic surface
\[
\{x \in {\mathbb R}^3 \: :\: f_3(x) = 0\}
\]
which is the algebraic closure of the boundary of
the pillow spectrahedron (the inner convex region)
\[
X=\{x \in {\mathbb R}^3 \: :\: f_2(x) \geq 0, \: f_3(x) \geq 0\}.
\]
In other words, the polynomial which vanishes along the boundary
of $X$ also vanishes outside of $X$, along
the Cayley cubic surface.
\end{example}

\section{Spectrahedral shadows and lifted LMIs}

We have seen that a spectrahedron is a closed basic semialgebraic set.
Moreover, it is a convex set. All spectrahedra are convex closed basic semialgebraic,
so one may then wonder conversely whether all convex closed basic semialgebraic sets
are spectrahedra. The answer is negative, even though we do not explain
why in this document.

\begin{proposition}[The TV screen is not a spectrahedron]\label{tvscreen}
The planar convex closed basic semialgebraic set
\[
\{x \in {\mathbb R}^2 \: :\: 1-x^4_1-x^4_2 \geq 0\}
\]
is not a spectrahedron.
\end{proposition}

Consequently, in order to represent convex closed basic semialgebraic sets, we have
to go beyond affine sections of the semidefinite cone. This motivates
the following definitions.

\begin{definition}[Lifted LMI, liftings]
A lifted LMI is a constraint
\[
F_0 + \sum_{k=1}^n x_k F_k + \sum_{l=1}^p u_l G_l \geq 0
\]
on a vector $x \in {\mathbb R}^n$, which implies additional variables
$u \in {\mathbb R}^p$ called liftings, and where matrices
$F_k \in {\mathbb S}^m$, $k=0,1,\ldots,n$
and $G_l \in {\mathbb S}^m$, $l=1,\ldots,p$ are given.
\end{definition}

\begin{definition}[Spectrahedral shadow, or lifted LMI set]
A spectrahedral shadow is the affine projection of a spectrahedron:
\[
\{x \in {\mathbb R}^n \: :\: F_0 + \sum_{k=1}^n x_k F_k + \sum_{l=1}^p u_l G_l \geq 0,
\: u \in {\mathbb R}^p \}
\]
where matrices $F_k \in {\mathbb S}^m$, $k=0,1,\ldots,n$ 
and $G_l \in {\mathbb S}^m$, $l=1,\ldots,p$ are given. 
\end{definition}

Spectrahedral shadows are also called semidefinite representable sets.

\begin{example}[The TV-screen is a spectrahedral shadow]
The planar convex closed basic semialgebraic set
\[
\{x \in {\mathbb R}^2 \: :\: 1-x^4_1-x^4_2 \geq 0\}
\]
can be expressed as the spectrahedral shadow
\[
\left\{x \in {\mathbb R}^2  \: :\: 
\left(\begin{array}{cc}1-u_1 & u_2 \\ u_2 & 1+u_1\end{array}\right) \geq 0, \quad
\left(\begin{array}{cc}1 & x_1 \\ x_1 & u_1\end{array}\right) \geq 0, \quad
\left(\begin{array}{cc}1 & x_2 \\ x_2 & u_2\end{array}\right) \geq 0,
\quad u \in {\mathbb R}^2\right\}.
\]
\end{example}

\begin{proposition}[Planar semialgebraic sets are spectrahedral shadows]\label{planar}
Every planar convex closed semialgebraic set is a spectrahedral shadow.
\end{proposition}

\begin{conjecture}[Are convex semialgebraic sets spectrahedral shadows ?]\label{convexalg}
Every convex closed semialgebraic set is a spectrahedral shadow.
\end{conjecture}

\section{SDP duality}

In this section we sketch some algorithms for linear semidefinite
programming (SDP). For notational simplicity we consider only the case
of a single semidefinite cone, with given data ${\mathcal A} : {\mathbb S}^n \to {\mathbb R}^m$,
$b \in {\mathbb R}^m$ and $C \in {\mathbb S}^n$.
We want to solve the primal SDP problem (\ref{psdp}) and
its dual SDP problem (\ref{dsdp}) reproduced here for the reader's convenience:
\[
\begin{array}{rcll}
p^* & = & \inf & \langle C,X \rangle \\
&& \mathrm{s.t.} & {\mathcal A}X = b \\
&&& X \geq 0
\end{array} \quad\quad
\begin{array}{rcll}
d^* & = & \sup & \langle b,y \rangle \\
&& \mathrm{s.t.} & Z=C-{\mathcal A}'y \\
&&& Z \geq 0.
\end{array}
\]
Define the feasibility sets
\[
P := \{X \in {\mathbb S}^n \: : \: {\mathcal A}X=b,\: X \geq 0\}, \quad 
D := \{y \in {\mathbb R}^m \: : \: C-{\mathcal A}'y \geq 0\}.
\]
Most of the algorithms for solving SDP problems make use of the
following elementary duality properties.

\begin{proposition}[Weak duality]\label{weakdual}
If $P$ and $D$ are nonempty, it holds $p^* \geq d^*$.
\end{proposition}

Indeed, if $P$ and $D$ are nonempty, there exist $X \in P$ and $y \in D$.
Letting $Z:=C-{\mathcal A}'y$ it holds
\[
p^*-d^* = \langle C,X \rangle - \langle b,y \rangle = \langle X,Z \rangle \geq 0
\]
since $X \geq 0$ and $Z \geq 0$.

\begin{proposition}[Strong duality]\label{strongdual}
If $P$ has nonempty interior and $D$ is nonempty, 
then the supremum $d^*$ is attained and $p^*=d^*$.
Similarly, if $D$ has nonempty interior and $P$ is nonempty,
then the infimum $p^*$ is attained and $p^*=d^*$.
\end{proposition}

Note in passing that whenever $X \geq 0$, $Z \geq 0$,
the scalar condition $p^*-d^*=\langle X,Z \rangle = 0$ is equivalent
to the nonsymmetric matrix condition $XZ = 0$ or to the
symmetric matrix condition $XZ+ZX=0$.

\section{Numerical SDP solvers}

Here we briefly describe numerical methods, implemented
in floating-point arithmetic, to solve SDP problems.
The most successful algorithms are primal-dual interior-point methods.

A triple $(X,y,Z)$ solves the primal-dual SDP problems (\ref{psdp})-(\ref{dsdp})
if and only if
\[
\begin{array}{c@{\quad}c}
{\mathcal A}X = b, \: X \geq 0 & \text{(primal feasibility)} \\
{\mathcal A}'y+Z = C, \: Z \geq 0 & \text{(dual feasibility)} \\
XZ = 0 & \text{(complementarity).}
\end{array}
\]
The key idea behind primal-dual interior-point methods is then
to consider the nonlinear system of equations
\begin{equation}\label{path}
\begin{array}{rcl}
{\mathcal A}X & = & b \\
{\mathcal A}'y + Z & = & C \\
XZ & = & \mu I_n
\end{array}
\end{equation}
parametrized in the scalar $\mu > 0$. These are
necessary and sufficiently optimality conditions for the strictly
convex problem
\[
\begin{array}{ll}
\inf & \langle C,X \rangle + \mu f(X) \\
\mathrm{s.t.} & {\mathcal A}X = b
\end{array}
\]
where $f(X):=-\log\det X$ is a barrier function for the semidefinite cone
(strictly convex and finite in its interior, and infinite elsewhere).
For this reason, scalar $\mu$ is called the barrier parameter.

If $P$ and $D$ have nonempty interior, it can be shown that 
for a given $\mu>0$, system (\ref{path}) has
a unique solution such that $X> 0$ and $Z > 0$,
and hence the set $\{(X(\mu),y(\mu),Z(\mu)) \: :\: \mu > 0\}$
defines a smooth curve parametrized by $\mu$, called the central path.
The interior-point algorithm consists then in applying Newton's method
to minimize a weighted sum of the linear objective function and
the barrier function, by following the central path,
letting $\mu\to 0$. Given initial feasible solutions, this generates a sequence of
feasible solutions such that the duality gap $\langle C,X \rangle
- \langle b,y \rangle$ is less that a given threshold $\epsilon>0$
after $O(\sqrt{n}\log\epsilon^{-1})$ iterations.  
Each Newton iteration requires:
\begin{itemize}
\item $O(n^2m)$ operations to evaluate the barrier function;
\item $O(n^3m)$ operations to evaluate and store its gradient;
\item $O(n^2m^2)$ operations to evaluate and store its Hessian;
\item $O(m^3)$ operations to solve the Newton linear system of equations.
\end{itemize}
A symmetric matrix of size $n$ has $O(n^2)$ independent entries, so
in general we may assume that $m=O(n^2)$ and hence the
dominating term in this rough complexity analysis
comes from the evaluation and storage
of the Hessian of the barrier function. Data sparsity
and block structure must be
exploited as much as possible in these steps.
It follows that a global worst-case asymptotic complexity
estimate for solving a dense SDP problem is $O(n^{6.5}\log\epsilon^{-1})$.
In practice the observed computational burden is much smaller,
but it strongly depends on the specific implementation
and on the problem structure.

Newton's method needs an initial feasible point, and if no such
point is available, an auxilliary SDP problem must be solved first.
An elegant approach to bypass the search of an initial point
consists of embedding the primal-dual
problem in a larger problem which is its own dual and
for which a trivial feasible starting point is known: this is
the so-called homogeneous self-dual embedding. A drawback of
this approach is that iterates are primal and dual feasible
for the original SDP problems only when the barrier parameter
vanishes.

The most successful semidefinite programming solvers are
implementations of primal-dual interior-point algorithms:
\begin{itemize}
\item {\tt SeDuMi}, {\tt SDPT3}, {\tt MOSEK}: homogenous self-dual embedding;
\item {\tt CSDP}, {\tt SDPA}: path-following predictor-corrector;
\item {\tt DSDP}: path-following with dual-scaling;
\end{itemize}
but there are also other implementations based on different
algorithms:
\begin{itemize}
\item {\tt LMILAB}: projective method;
\item {\tt PENSDP}: penalty and augmented Lagrangian.
\end{itemize}
There exist parallel implementations of {\tt CSDP} and {\tt SDPA}.

Most of these solvers are available under {\tt Matlab}, and they
are interfaced through the parsers {\tt YALMIP} and {\tt cvx}. Some
elementary SDP solver is available under {\tt Scilab} and
{\tt Sage}, and {\tt cvxopt} is a {\tt Python} interface
with some SDP features. The solver {\tt CSDP} can
be embedded in C language, the solver {\tt SDPA}
is also available with {\tt Python} interface,
and {\tt PENSDP} is available as a standalone
solver or called in {\tt Fortran} or C language.

\section{Rigorous SDP solvers}

The numerical methods described in the previous sections
are implementable in floating-point arthimetic, but
very little is known about backward stability of
these algorithms. More annoyingly, it is difficult
to estimate or bound the conditioning of a SDP problem,
which implies that none of these numerical solvers
can provide a priori guarantees about the quality of their output,
even for a restricted problem class.

To address this issue, various strategies can be followed:
\begin{itemize}
\item multiprecision arithmetic;
\item interval arithmetic;
\item symbolic computation.
\end{itemize}

Higher precision or arbitrary precision arithmetic allows to deal
with better floating-point approximations of real numbers,
at the price of an increased computational burden. Currently,
the solver {\tt SDPA} is available in quad-double precision,
double-double precision and arbitrary precision
arithmetic.

Interval arithmetic can be used to obtain rigorous bounds
on the output of numerical SDP solvers. A {\tt Matlab}
implementation of a verified SDP solver is {\tt VSDP}.
It relies on the {\tt Intlab} toolbox for interval
computations.

Symbolic computation can be used to solve SDP problems
exactly, by solving (e.g. with Gr\"obner basis techniques)
the quadratic system of equations arising from optimality
conditions. Alternatively, feasible points in spectrahedra
can be obtained by techniques for finding real solutions
of systems of polynomial equations and inequalities.

To justify further the need for these techniques,
note first that there are SDP problems with integer data
with no solution among the rationals:
\begin{example}[Irrational optimal solution]\label{irrat1}
The problem
\[
\begin{array}{rl}
\sup & y \\
\mathrm{s.t.} & \left(\begin{array}{cc} 1 & y \\ y & 2 \end{array}\right) \geq 0
\end{array}
\]
has solution $y^*=\sqrt{2}$.
\end{example}
\begin{example}[Irrational spectrahedron]\label{irrat2}
\[
\{y \in {\mathbb R} \: :\: \left(\begin{array}{cc}1&y\\y&2\end{array}\right)
\geq 0, \: \left(\begin{array}{cc}2y&2\\2&y\end{array}\right) \geq 0\}
\:=\: \{\sqrt{2}\}.
\]
\end{example}

In general, exact solutions of SDP problems must be found
in algebraic extensions of the ground field of the input data.
Recall that when both primal and dual problems
have nonempty interiors, solutions $(X,y)$ are characterized
by the optimality conditions (\ref{path}) with $\mu=0$, i.e.
\begin{equation}\label{optimal}
\begin{array}{rcl}
\langle A_k,X \rangle & = & b_k, \quad k=1,\ldots,m \\
X(C-\sum_{k=1}^m y_k A_k) & = & 0_n.
\end{array}
\end{equation}
This is a system of $m+n(n+1)/2$ real linear and quadratic
equations in $m+n(n+1)/2$ real variables. If we have a basis for
the nullspace of the operator $\mathcal A$, we can remove the first $m$
equality constraints and derive a system of $n(n+1)/2$ quadratic
equations in $n(n+1)/2$ variables.

\begin{example}[Irrational optimal solution, again]\label{irrat1b}
Optimality conditions for the problem of Example \ref{irrat1} are as follows:
\[
\begin{array}{rcl}
-x_{21} & = & 1 \\
x_{11}+yx_{21} & = & 0 \\
x_{21}+yx_{22} & = & 0 \\
yx_{21}+2x_{22} & = & 0 
\end{array}
\Longleftrightarrow
\begin{array}{rcl}
x_{11} & = & \pm\sqrt{2} \\
x_{21} & = & -1 \\
x_{22} & = & \pm\frac{\sqrt{2}}{2} \\
y & = & \pm\sqrt{2}
\end{array}
\]
from which it follows that the primal-dual optimal solution is
\[
X^* = \left(\begin{array}{cc}
\sqrt{2} & -1 \\ -1 & \frac{\sqrt{2}}{2}
\end{array}\right) \quad
y^* = \sqrt{2}.
\]
\end{example}

In the classical Turing machine model of computation,
an integer number $N$ is encoded in binary notation,
so that its bit size is $\log_2 N+1$. The following
spectrahedron with integer coefficients has points
with exponential bit size:
\begin{example}[Exponential spectrahedron]\label{expo}
Any point in the spectrahedron
\[
\{y \in {\mathbb R}^m \: :\: \left(\begin{array}{cc}
1 & 2 \\ 2 & y_1 \end{array}\right) \geq 0, \:
\left(\begin{array}{cc}
1 & y_1 \\ y_1 & y_2
\end{array}\right) \geq 0, \:\cdots\:,
\left(\begin{array}{cc}
1 & y_{m-1} \\ y_{m-1} & y_m
\end{array}\right) \geq 0\}
\]
satisfies $y_m \geq 2^{2^m}$.
\end{example}

\begin{example}[Algebraic solution]\label{algebraic}
Consider the problem 
\[
\begin{array}{ll}
\sup & y_1 + y_2 + y_3 \\
\mathrm{s.t.} & \left(\begin{array}{cccc}
1+y_3 & y_1+y_2 & y_2 & y_2+y_3\\
y_1+y_2 & 1-y_1 & y_2-y_3 & y_2 \\
y_2 & y_2-y_3 & 1+y_2 & y_1+y_3 \\
y_2+y_3 & y_2 & y_1+y_3 & 1-y_3
\end{array}\right) \geq 0.
\end{array}
\]
Optimality conditions (\ref{optimal})
yield 13 equations in 13 unknowns. Using Gr\"obner basis techniques,
it is found that these equations have 26 complex solutions.
The optimal first variable $y^*_1$ is the root of a degree 26
univariate polynomial with integer coefficients. This polynomial
factors into a degree 16 term
\[
403538653715069011 y^{16}_1 - 2480774864948860304 y^{15}_1
+ \cdots + 149571632340416
\]
and a degree 10 term
\[
2018 y^{10}_1 - 12156 y^9_1 + 17811 y^8_1 + \cdots - 163
\]
both irreducible in ${\mathbb Q}[y_1]$.
The optimal solution $y^*_1$ is therefore an algebraic number
of degree 16 over ${\mathbb Q}$, and it can be checked that
it is also the case for the other 12 optimal coordinates
$y^*_2,y^*_3,x^*_{11},x^*_{21},\ldots,x^*_{44}$.
\end{example}

The above examples indicate that it can be quite costly
to solve an SDP problem exactly. The algebraic degree of
an SDP problem is the degree of the algebraic extension
of the problem data coefficient field
over which the solutions should be found. Even for
small $n$ and $m$, this number can be very large.

\section{Notes and references}

References on convex analysis
are \cite{r70} and \cite{hl01}.
See \cite{bv05} for an elementary introduction
to convex optimization, and \cite{bn01}
for a more advanced treatment aimed at
applied mathematicians and engineers.
Systems control applications of linear matrix inequalities 
are described in \cite{befb94}.
Good historical surveys on SDP
are \cite{vb96} and \cite{t01}. Classifications
of sets and functions that can be represented by
affine sections and projections of LP, QP and SDP
cones can be found in \cite{nn94}, \cite{bn01}
and \cite{n06}. Elementary concepts of algebraic
geometry (algebraic sets, semialgebraic sets)
are surveyed in \cite{clo07}, and connections
between SDP, convex geometry and algebraic geometry
are explored in \cite{bpt13}.
Proposition \ref{csa} is proved in \cite[Theorem 20]{r06}.
Example \ref{ex:pillow} comes from the SDP relaxation
of a 3-dimensional MAXCUT problem, a classical
problem of combinatorial optimization, see \cite{lr05}
and also \cite[Example 2]{nrs10} for the link
with the Cayley cubic.
The example of Proposition \ref{tvscreen} was studied
in \cite{hv07}. The proof of Proposition \ref{planar}
can be found in \cite{s12}. Conjecture \ref{convexalg},
a follow-up of a question posed in \cite[Section 4.3.1]{n06},
can be found in \cite{hn09}.
A basic account of semidefinite
programming duality (Propositions \ref{strongdual}
and \ref{weakdual}), as well as Examples \ref{irrat2} and \ref{expo}
can be found in \cite[Section 2]{lr05}.
Techniques of real algebraic geometry for finding
rational points in convex semialgebraic sets are
described in \cite{sz10}.
Example \ref{algebraic} is taken from \cite[Example 4]{nrs10},
which describes an approach to quantifying the complexity of solving
exactly an SDP problem. Hyperlinks to SDP solvers can be found easily,
and the online documentation of the interface {\tt YALMIP}
contains many pointers to webpages and software packages.


\chapter{Finite-dimensional polynomial optimization}\label{chap:polyopt}

\section{Measures and moments}

Let $X$ be a compact subset of the Euclidean space ${\mathbb R}^n$.
Let ${\mathscr B}(X)$ denotes
the Borel $\sigma$-algebra, defined as the smallest collection of subsets
of $X$ which contains all open sets.

\begin{definition}[Signed measure]
A signed measure is a function $\mu : {\mathscr B}(X) \to {\mathbb R} \cup \{\infty\}$
such that $\mu(\emptyset)=0$ and $\mu(\cup_{k\in{\mathbb N}} X_k) = 
\sum_{k\in{\mathbb N}}\mu(X_k)$ for any pairwise disjoint $X_k \in {\mathscr B}(X)$.
\end{definition}

\begin{definition}[Positive measure]
A positive measure is a signed measure which takes only nonnegative values.
\end{definition}

Positive measures on the Borel $\sigma$-algebra
are often called Borel measures, and positive measures which take
finite values on compact sets are often called Radon measures.

\begin{definition}[Support]
Given a measure $\mu$ its support $\mathrm{spt}\:\mu$
is the closed set of all points $x$ such that $\mu(A)\neq0$
for every neighborhood $A$ of $x$. We say that $\mu$ is
supported on a set $A$ whenever $\mathrm{spt}\:\mu \subset A$.
\end{definition}

\begin{definition}[Probability measure]
A probability measure $\mu$ on $X$ is a positive measure such that $\mu(X)=1$.
\end{definition}

Let us denote by ${\mathscr M}_+(X)$ the cone
of positive measures supported on $X$, and 
by ${\mathscr P}(X)$ the set
of probability measures supported on $X$. Geometrically,
${\mathscr P}(X)$ is an affine section of ${\mathscr M}_+(X)$.

\begin{example}[Lebesgue measure]
The Lebesgue measure on ${\mathbb R}^n$, also called uniform measure,
denoted $\lambda$, is a positive measure returning the volume of a set $A$.
For instance, when $n=1$ and $a\leq b$, $\lambda([a,b]) = b-a$.
\end{example}

\begin{example}[Dirac measure]
The Dirac measure at $x=\xi$, denoted $\_\xi(dx)$ or $\delta_{x=\xi}$,
is a probability measure such that $\delta_\xi(A)=1$ if $\xi \in A$,
and $\delta_\xi(A)=0$ if $\xi \notin A$.
\end{example}

For a given compact set $X \subset {\mathbb R}^n$, let ${\mathscr M}(X)$ denote the Banach space of
signed measures supported on $X$, so that a measure $\mu \in {\mathscr M}(X)$
can be interpreted as a function that
takes any subset of $X$ and returns a real number. Alternatively, elements of
${\mathscr M}(X)$ can be interpreted as linear functionals acting on the Banach space of
continuous functions ${\mathscr C}(X)$, that is, as elements of the
dual space ${\mathscr C}(X)'$, see Definition \ref{dualspace}.
The action of a measure $\mu \in {\mathscr M}(X)$
on a test function $v \in {\mathscr C}(X)$
can be modeled with the duality pairing
\[
\langle v,\:\mu \rangle \: := \int_X v(x)\,d\mu(x).
\]
Let us denote by ${\mathscr C}_+(X)$ the cone of positive continuous functions on $X$,
whose dual can be identified to the cone of positive measures on $X$, i.e.
${\mathscr M}_+(X)={\mathscr C}_+(X)'$.

\begin{definition}[Monomial]
Given a real vector $x \in {\mathbb R}^n$ and an integer vector $\alpha \in {\mathbb N}^n$,
a monomial is defined as
\[
x^{\alpha} := \prod_{k=1}^n x^{\alpha_k}_k.
\]
\end{definition}

The degree of a monomial with exponent $\alpha \in {\mathbb N}^n$ is
equal to $|\alpha|:=\sum_{k=1}^n \alpha_k$.

\begin{definition}[Moment]
Given a measure $\mu \in {\mathscr M}(X)$, the real number
\begin{equation}\label{mom}
y_{\alpha} := \int_X x^{\alpha} \mu(dx)
\end{equation}
is called its moment of order $\alpha \in {\mathbb N}^n$.
\end{definition}

\begin{example}
For $x \in {\mathbb R}^2$, second order moments are
\[
y_{20}=\int x^2_1\mu(dx), \:\: y_{11}=\int x_1x_2\mu(dx), \:\: y_{02}=\int x^2_2\mu(dx).
\]
\end{example}

The sequence $(y_{\alpha})_{\alpha \in {\mathbb N}^n}$ is called the
sequence of moments of the measure $\mu$, and given $d \in {\mathbb N}$,
the truncated sequence $(y_{\alpha})_{|\alpha|\leq d}$ is the
vector of moments of degree $d$. 

\begin{definition}[Representing measure]
If $y$ is the sequence of moments of a measure $\mu$, i.e. if
identity (\ref{mom}) holds for all $\alpha \in {\mathbb N}^n$,
we say that $\mu$ is a representing measure for $y$.
\end{definition}

A basic problem in
the theory of moments concerns the characterization of (infinite
or truncated) sequences that are moments of some measure.
Practically speaking,
instead of manipulating a measure, which is a rather abstract
object, we manipulate its moments. Indeed, a measure on a
compact set is uniquely determined by the (infinite) sequence
of its moments.

\section{Riesz functional, moment and localizing matrices}

Measures on $X \subset {\mathbb R}^n$ are manipulated
with their moments, i.e. via their actions on monomials.
The choice of monomials $(x^{\alpha})_{\alpha}$ is motivated
mainly for notational and simplicity reasons. In particular,
the product of two monomials is a monomial, i.e.
$x^{\alpha}x^{\beta}=x^{\alpha+\beta}$. Any other choice
of basis $(b_{\alpha}(x))_{\alpha}$ would be appropriate
to manipulate measures,
as soon as the basis is dense w.r.t. the supremum norm
in the space of continuous functions ${\mathscr C}(X)$. Numerically
speaking, other bases than monomials may be more appropriate,
but we do not elaborate further on this issue in this document.

In order to manipulate functions in ${\mathscr C}(X)$ we use polynomials.
A polynomial $p \in {\mathbb R}[x]$ of degree $d \in {\mathbb N}$
is understood as a linear combination of monomials:
\[
p(x) := \sum_{|\alpha|\leq d} p_{\alpha} x^{\alpha}
\]
and $p:=(p_{\alpha})_{|\alpha|\leq d}$ is the vector
of its coefficients in the monomial basis $(x^{\alpha})_{\alpha}$.
Note that we use the same notation for a polynomial and for its
vector of coefficients when no ambiguity is possible.
Otherwise we use the notation $p(x)$ to emphasize that we
deal the polynomial as a function, not as a vector.

\begin{example}\label{ex:poly}
The polynomial
\[
x \in {\mathbb R}^2 \:\:\mapsto\:\: p(x)=1+2x_2+3x^2_1+4x_1x_2
\]
has a vector of coefficients $p \in {\mathbb R}^6$ with entries
$p_{00}=1$, $p_{10}=0$, $p_{01}=2$, $p_{20}=3$, $p_{11}=4$, $p_{02}=0$.
\end{example}

\begin{definition}[Riesz functional]
Given a sequence $y = (y_{\alpha})_{\alpha \in {\mathbb N}^n}$,
we define the Riesz linear functional $\ell_y : {\mathbb R}[x] \to {\mathbb R}$
which acts on polynomials $p(x) = \sum_{\alpha} p_{\alpha} x^{\alpha}$
as follows: $\ell_y(p(x)) = \sum_{\alpha} p_{\alpha} y_{\alpha}$.
\end{definition}

We can interpret the Riesz functional as an operator
that linearizes polynomials. If sequence $y$ has a representing measure $\mu$,
integration of a polynomial $p$ w.r.t. $\mu$ is obtained by applying the
Riesz functional $\ell_y$ on $p$, since
\[
\ell_y(p) = \sum_{\alpha} p_{\alpha} y_{\alpha} =
\sum_{\alpha} p_{\alpha} \int x^{\alpha} \mu(dx) =
\int \sum_{\alpha} p_{\alpha} x^{\alpha} \mu(dx) =
\int p(x)\mu(dx).
\]
Note that formally the Riesz functional is the linear form
$p(x) \mapsto \ell_y(p(x))$ and its existence is independent
of the choice of basis to represent polynomial $p(x)$.
However, for notational simplicity, we use the monomial basis
and hence we represent explicitly the Riesz functional with
the inner product of the vector $(p_{\alpha})_{|\alpha|\leq d}$
of coefficients of the polynomial with the truncated sequence
$(y_{\alpha})_{|\alpha|\leq d}$.

\begin{example}
For the polynomial of Example \ref{ex:poly}, the Riesz functional
reads
\[
p(x)=1+2x_2+3x^2_1+4x_1x_2 \:\:\mapsto\:\: \ell_y(p) = y_{00}+2y_{01}+3y_{20}+4y_{11}.
\]
\end{example}

If we apply the Riesz functional on the square of a polynomial $p(x)$,
then we obtain a form which is quadratic in the coefficients of $p(x)$:
\begin{definition}[Moment matrix]
The moment matrix of order $d$ is the Gram matrix of
the quadratic form $p(x) \mapsto \ell_y(p^2(x))$ where
polynomial $p(x)$ has degree $d$, i.e.
the matrix $M_d(y)$
such that $\ell_y(p^2(x))=p'M_d(y)p$.
\end{definition}

\begin{example}
If $n=2$ then
\[
M_0(y) = y_{00}, \quad
M_1(y) = \left(\begin{array}{ccc}
y_{00} & y_{10} & y_{01} \\
y_{10} & y_{20} & y_{11} \\
y_{01} & y_{11} & y_{02}
\end{array}\right), \quad
M_2(y) = \left(\begin{array}{cccccc}
y_{00} & y_{10} & y_{01} & y_{20} & y_{11} & y_{02} \\
y_{10} & y_{20} & y_{11} & y_{30} & y_{21} & y_{12} \\
y_{01} & y_{11} & y_{02} & y_{21} & y_{12} & y_{03} \\
y_{20} & y_{30} & y_{21} & y_{40} & y_{31} & y_{22} \\
y_{11} & y_{21} & y_{12} & y_{31} & y_{22} & y_{13} \\
y_{02} & y_{12} & y_{03} & y_{22} & y_{13} & y_{04}
\end{array}\right).
\]
\end{example}

Note that $M_d(y) \in {\mathbb S}^{n+d \choose n}$
where 
\[
{n+d \choose n} = {n+d \choose d} = \frac{(n+d)!}{n!\:d!}
\]
is the number of monomials of $n$ variables of degree at most $d$.
The rows and columns of the moment matrix 
are indexed by vectors $\alpha \in {\mathbb N}^n$ and $\beta \in {\mathbb N}^n$.
Inspection reveals that indeed the entry $(\alpha,\beta)$ in the moment matrix
is the moment $y_{\alpha+\beta}$. 
By construction, the moment matrix $M_d(y)$ is symmetric and linear in $y$.

If we apply the Riesz functional on the product of the square of a polynomial $p(x)$
of degree $d$ with a given polynomial $q(x)$, then we obtain a form which is quadratic
in the coefficients of $p(x)$.

\begin{definition}[Localizing matrix]
Given a polynomial $q(x)$, its localizing matrix of order $d$ is the Gram matrix of
the quadratic form $p(x) \mapsto \ell_y(q(x)p^2(x))$ where polynomial $p(x)$
has degree $d$, i.e. the matrix $M_d(q\:y)$
such that $\ell_y(q(x)p^2(x))=p'M_d(q\:y)p$.
\end{definition}

Note that we use the notation $M_d(q\:y)$ to emphasize the fact that the
localizing matrix is bilinear in $q$ and $y$. When polynomial
$q(x)=\sum_\alpha q_{\alpha} x^{\alpha}$ is given,
matrix $M_d(q\:y)$ is symmetric and linear in $y$.
The localizing matrix can be interpreted as a linear combination
of moment matrices, in the sense that its entry $(\alpha,\beta)$
is equal to $\sum_{\gamma} q_{\gamma} y_{\alpha+\beta+\gamma}$.

\begin{example}
If $n=2$ and $q(x)=1+2x_1+3x_2$ then
\[
M_1(q\:y)=\left(\begin{array}{ccc}
y_{00}+2y_{10}+3y_{01} & y_{10}+2y_{20}+3y_{11} &  y_{01}+2y_{11}+3y_{02} \\
y_{10}+2y_{20}+3y_{11} & y_{20}+2y_{30}+3y_{21} & y_{11}+2y_{21}+3y_{12} \\
y_{01}+2y_{11}+3y_{02} & y_{11}+2y_{21}+3y_{12} & y_{02}+2y_{12}+3y_{03}
\end{array}\right).
\]
\end{example}

Finally, given an infinite-dimensional sequence $y$, 
let us denote the infinite dimensional moment and localized matrices, or linear
operators, as follows
\[
M(y):=M_{\infty}(y), \quad M(q\:y):=M_{\infty}(q\:y).
\]

\section{Linking measures and moments}\label{sec:repres}

The matrices just introduced allow to
explicitly model the constraint that a sequence $y$ has a representing 
measure $\mu$ on a compact basic semialgebraic set $X$.
Under a mild assumption on the representation of $X$,
it turns out that this constraint is an infinite-dimensional LMI.

\begin{assumption}[Compactness]\label{compact}
Assume that $X$ is a compact basic semialgebraic set
\[
X:=\{x \in {\mathbb R}^n \: :\: p_k(x) \geq 0, \:k=1,\ldots,n_X\}
\]
for given $p_k \in {\mathbb R}[x]$, $k=1,\ldots,n_X$.
Moreover, assume that one of the polynomial inequalities $p_k(x)\geq 0$ is of the
form $R-\sum_{i=1}^n x^2_i\geq 0$ where $R$ is a sufficiently large positive constant.
\end{assumption}

On the one hand, Assumption \ref{compact} is a little bit stronger than
compactness of $X$. On the other hand, if we assume only that $X$ is compact,
this is without loss of generality that a constraint 
can be added to the description of $X$ so that
Assumption \ref{compact} is satisfied.

\begin{proposition}[Putinar's Theorem]\label{putinar}
Let set $X$ satisfy Assumption \ref{compact}. Then
sequence $y$ has a representing measure
in ${\mathscr M}_+(X)$ if and only if $M(y)\geq 0$ and $M(p_k\:y) \geq 0$, $k=1,\ldots,n_X$.
\end{proposition}

Note that if we have an equality constraint $p_k(x)=0$ instead of
an inequality constraint in the definition of $X$,
the corresponding localizing constraint becomes $M(p_k\:y)=0$, which
is a set of linear equations in $y$.

Since matrices $M(y)$ and $M(p_k\:y)$ are symmetric and linear in $y$,
sequences with representing measures belong to an infinite-dimensional
spectrahedron, following the terminology introduced in Section \ref{lmi}.
To manipulate these objects, we will consider finite-dimensional
truncations.

\section{Measure LP}

Let $X$ be the compact basic semialgebraic set defined above
for given polynomials $p_k \in {\mathbb R}[x]$, $k=1,\ldots,n_X$,
and satisfying Assumption \ref{compact}.
Let $p_0 \in {\mathbb R}[x]$ be a given a polynomial.
Consider the optimization problem consisting of
minimizing $p_0$ over $X$, namely
\begin{equation}\label{poly}
\begin{array}{rcll}
p^* & = & \min & p_0(x) \\
&& \mathrm{s.t.} & p_k(x)\geq 0, \: k=1,\ldots,n_K.
\end{array}
\end{equation}
The above minimum is w.r.t. $x \in {\mathbb R}^n$ and since
we assume that $X$ is compact, the minimum is attained
at a given point $x^* \in X$.

We do not have any convexity property on $p_0$ or $X$, so that
problem (\ref{poly}) may feature several local minima, and
possibly several global minima. In the sequel we describe
a hierarchy of LMI relaxations of increasing size, indexed by
a relaxation order, and that generates an asymptotically convergent
mononotically nondecreasing sequence of lower bounds on $p^*$.

The key idea is to notice
that nonconvex polynomial optimization problem (\ref{poly})
over (finite-dimensional set) $X \subset {\mathbb R}^n$
is equivalent to a linear, hence convex, optimization problem
over the (infinite-dimensional set) of probability measures
supported on $X$. More specifically, consider the problem
\begin{equation}\label{linopt}
\begin{array}{rcll}
p^*_M & = & \inf & \int p_0(x)\mu(dx) \\
&& \mathrm{s.t.} & \mu(X)=1\\
&&& \mu \in {\mathscr M}_+(X)
\end{array}
\end{equation}
which is linear in the decision variable $\mu$,
a probability measure supported on $X$.

\begin{proposition}[Measure LP formulation of polynomial optimization]
The infimum in LP problem (\ref{linopt}) is attained,
and $p^*_M=p^*$.
\end{proposition}

The proof is immediate: for any feasible $\xi \in X$, it holds
$p_0(\xi)=\int p_0(x)\mu(dx)$ for the Dirac measure $\mu=\delta_{\xi}$,
showing $p^*\geq p^*_M$. Conversely, as
$p_0(x)\geq p^*$ for all $x \in X$, it holds $\int_X p_0(x)\mu(dx) \geq \int_X p^*\mu(dx)
= p^* \int_X \mu(dx)= p^*$ since $\mu$ is a probability measure,
which shows that $p^*_M\geq p^*$. It follows that $p^*_M=p^*$ and that 
the infimum in problem (\ref{linopt}) is attained
by a Dirac measure $\mu=\delta_{x^*}$ where $x^*$ is a global
optimum of problem (\ref{poly}).

\section{Moment LP}\label{sec:momlp}

In Section \ref{chap:finlp} we studied LP problems (\ref{plp})
in finite-dimensional cones. In the context of polynomial optimization,
we came up with infinite-dimensional LP (\ref{linopt}) which is a
special instance of the measure LP
\begin{equation}\label{infp}
\begin{array}{rcll}
p^* & = & \inf & \langle c,\mu \rangle \\
&& \mathrm{s.t.} & {\mathcal A}\mu = b \\
&&& \mu \in {\mathscr M}_+(X)
\end{array}
\end{equation}
where the decision variable $\mu$ is in the cone
of nonnegative measures supported on $X$, a given compact
subset of ${\mathbb R}^n$. Linear operator ${\mathcal A} : {\mathscr M}(X) \to {\mathbb R}^m$
takes a measure and returns an $m$-dimensional vector
of real numbers. Vector $b \in {\mathbb R}^m$ is given.
The objective function is the duality pairing between
a given continuous function $c \in {\mathscr C}(X)$ and $\mu$.
Problem (\ref{infp}) has a dual (or more rigorously, a predual) problem
in the cone of nonnegative functions, but we will not describe it in
this document.

If the linear operator $\mathcal A$ is described
through given continuous functions $a_j \in {\mathscr C}(X)$,
$j=1,\ldots,m$ we can write the LP problem (\ref{infp})
more explicitly as
\begin{equation}\label{pmea}
\begin{array}{rcll}
p^* & = & \inf & \int_X c(x)\mu(dx) \\
&& \mathrm{s.t.} & \int_X a_j(x)\mu(dx) = b_j, \:\: j=1,\ldots,m \\
&&& \mu \in {\mathscr M}_+(X).
\end{array}
\end{equation}
Now suppose that all the functions are
polynomials, i.e. $a_j(x) \in {\mathbb R}[x]$, $j=1,\ldots,m$,
$c(x) \in {\mathbb R}[x]$,
so that measure $\mu$ can be manipulated via the sequence
$y:=(y_{\alpha})_{\alpha \in {\mathbb N}}$ of its moments (\ref{mom}).
The measure LP (\ref{pmea}) becomes a moment LP
\begin{equation}\label{pmom}
\begin{array}{rcll}
p^* & = & \inf & \sum_{\alpha} c_{\alpha}y_{\alpha} \\
&& \mathrm{s.t.} & \sum_{\alpha} a_{j\alpha}y_{\alpha} = b_j, \:\: j=1,\ldots,m \\
&&& \text{$y$ has a representing measure $\mu \in {\mathscr M}_+(X)$}
\end{array}
\end{equation}
called a generalized problem of moments.

The idea is then to use the explicit LMI conditions of Section \ref{sec:repres}
to model the constraints that a sequence has a representing measure.
If the semialgebraic set
\[
X := \{x \in {\mathbb R}^n \: :\: p_k(x) \geq 0, \: k=1,\ldots,n_{X}\}
\]
satisfies Assumption \ref{compact}, problem (\ref{pmom}) becomes 
\[
\begin{array}{rcll}
p^* & = & \inf & \sum_{\alpha} c_{\alpha}y_{\alpha} \\
&& \mathrm{s.t.} & \sum_{\alpha} a_{j\alpha}y_{\alpha} = b_j, \:\: j=1,\ldots,m \\
&&& M(y) \geq 0, \: M(p_k\:y) \geq 0, \: k=1,\ldots,n_X
\end{array}
\]
where the constraints $\sum_{\alpha} a_{j\alpha}y_{\alpha} = b_j$, $j=1,\ldots,m$
model finitely many linear constraints on infinitely many decision variables.
In the sequel, we will consider finite-dimensional truncations of this problem,
and generate a hierarchy of LMI relaxations called Lasserre's hierarchy
in the context of polynomial optimization.

Moment LP (\ref{pmom}) has a dual in the cone of positive polynomials,
and finite-dimensional truncations of this problem correspond to
the search of polynomial sum-of-squares representations, which can
be formulated with a hierarchy of dual LMI problems, but we will
not elaborate more on this point in this document.

\section{Lasserre's LMI hierarchy}\label{sec:hierarchy}

Now remark that LP problem (\ref{linopt}) is a special
instance of the moment LP problem (\ref{pmea})
with data $c(x)=p_0(x)=\sum_{\alpha} {p_0}_{\alpha} x^{\alpha}$, $a(x)=1$, $b=1$,
so that, as in Section \ref{sec:momlp},
problem (\ref{poly}) can be equivalently written as
\[
\begin{array}{rcll}
p^* & = & \inf & \sum_{\alpha} {p_0}_{\alpha}y_{\alpha} \\
&& \mathrm{s.t.} & y_0 = 1 \\
&&& M(y) \geq 0, \: M(p_k\:y) \geq 0, \: k=1,\ldots,n_X.
\end{array}
\]
Let us denote by $r_k$ the smallest integer not less than half
the degree of polynomial $p_k$, $k=0,1,\ldots,n_X$, and let
$r_X:=\max\{1,r_1,\ldots,r_{n_X}\}$.
For $r \geq r_X$, consider {\bf Lasserre's LMI hierarchy}
\begin{equation}\label{hierarchy}
\begin{array}{rcll}
p^*_r & = & \inf & \sum_{\alpha} {p_0}_{\alpha}y_{\alpha} \\
&& \mathrm{s.t.} & y_0 = 1 \\
&&& M_r(y) \geq 0, \: M_{r-r_k}(p_k\:y) \geq 0, \: k=1,\ldots,n_X.
\end{array}
\end{equation}
The LMI constraints in this problem are truncated, or relaxed
versions of the infinite-dimensional LMI constraints
of Proposition \ref{putinar}. When the {\bf relaxation order} $r \in {\mathbb N}$
tends to infinity, we obtain the following result.

\begin{proposition}[Lasserre's LMI hierarchy converges]\label{lasserre}
It holds $p^*_r  \leq p^*_{r+1} \leq p^*$ and $\lim_{r\to\infty} p^*_r=p^*$.
\end{proposition}

Lasserre's LMI relaxations (\ref{hierarchy}) can be solved with
semidefinite programming, see Chapter \ref{chap:finlp},
and this provides us with a monotonically nondecreasing sequence of lower bounds on
the global minimum of nonconvex polynomial optimization problem (\ref{poly}).

\begin{proposition}[Generic finite convergence]\label{generic}
In the finite-dimensional space of coefficients
of polynomials $p_k$, $k=0,1,\ldots,n_X$ defining problem (\ref{poly}),
there is a low-dimensional algebraic set
which is such that if we choose an instance of problem (\ref{poly}) outside
of this set, then Lasserre's LMI relaxations have finite convergence,
i.e. there exists a finite $r^*$ such that $p^*_r=p^*$ for all $r\geq r^*$.
\end{proposition}

Equivalently, finite convergence occurs under arbitrary small perturbations
of the data of problem (\ref{poly}), and problems for which finite convergence
does not occur are exceptional and degenerate in some sense.

\begin{example}\label{ex:polyopt}
Consider the polynomial optimization problem
\[
\begin{array}{rcll}
p^* & = & \min & -x_2 \\
&& \mathrm{s.t.} & 3-2x_2-x^2_1-x^2_2 \geq 0 \\
&&& -x_1-x_2-x_1x_2 \geq 0 \\
&&& 1+x_1x_2 \geq 0
\end{array}
\]
where the minimum is w.r.t. $x \in {\mathbb R}^2$.
The first LMI relaxation is
\[
\begin{array}{rcll}
p^*_1 & = & \min & -y_{01} \\
&& \mathrm{s.t.} & y_{00} = 1 \\
&&& \left(\begin{array}{ccc}
y_{00}&y_{10}&y_{01}\\
y_{10}&y_{20}&y_{11}\\
y_{01}&y_{11}&y_{02}\end{array}\right) \geq 0 \\
&&& 3y_{00}-2y_{01}-y_{20}-y_{02} \geq 0 \\
&&& -y_{10}-y_{01}-y_{11} \geq 0 \\
&&& y_{00}+y_{11} \geq 0
\end{array}
\]
and the second LMI relaxation is
\[
\begin{array}{rcll}
p^*_2 & = & \min & -y_{01} \\
&& \mathrm{s.t.} & y_{00} = 1 \\
&&& \left(\begin{array}{cccccc}
y_{00}&y_{10}&y_{01}&y_{20}&y_{11}&y_{02}\\
y_{10}&y_{20}&y_{11}&y_{30}&y_{21}&y_{12}\\
y_{01}&y_{11}&y_{02}&y_{21}&y_{12}&y_{03}\\
y_{20}&y_{30}&y_{21}&y_{40}&y_{31}&y_{22}\\
y_{11}&y_{21}&y_{12}&y_{31}&y_{22}&y_{13}\\
y_{02}&y_{12}&y_{03}&y_{22}&y_{13}&y_{04}
\end{array}\right) \geq 0 \\ 
&&& \left(\begin{array}{ccc}
3y_{00}-2y_{01}-y_{20}-y_{02}&3y_{10}-2y_{11}-y_{30}-y_{12}&
3y_{01}-2y_{02}-y_{21}-y_{03}\\
3y_{10}-2y_{11}-y_{30}-y_{12}&3y_{20}-2y_{21}-y_{40}-y_{22}&
3y_{11}-2y_{12}-y_{31}-y_{13}\\
3y_{01}-2y_{02}-y_{21}-y_{03}&3y_{11}-2y_{12}-y_{31}-y_{13}&
3y_{02}-2y_{03}-y_{22}-y_{04}
\end{array}\right) \geq 0 \\
&&& \left(\begin{array}{ccc}
-y_{10}-y_{01}-y_{11}&-y_{20}-y_{11}-y_{21}&-y_{11}-y_{02}-y_{12}\\
-y_{20}-y_{11}-y_{21}&-y_{30}-y_{21}-y_{31}&-y_{21}-y_{12}-y_{22}\\
-y_{11}-y_{02}-y_{12}&-y_{21}-y_{12}-y_{22}&-y_{12}-y_{03}-y_{13}
\end{array}\right) \geq 0 \\
&&& \left(\begin{array}{ccc}
y_{00}+y_{11}&y_{10}+y_{21}&y_{01}+y_{12}\\
y_{10}+y_{21}&y_{20}+y_{31}&y_{11}+y_{22}\\
y_{01}+y_{12}&y_{11}+y_{22}&y_{02}+y_{13}
\end{array}\right) \geq 0. \\
\end{array}
\]
It can be checked that $p^*_1 = -2 \leq p^*_2 = p^* = -\frac{1+\sqrt{5}}{2}$.
Note that Assumption \ref{compact} is satisfied for this example,
since the constraint $3-2x_2-x^2_1-x^2_2 \geq 0$ certifies
boundedness of the feasibility set.
\end{example}

\section{Global optimum recovery}\label{recover}

From Proposition \ref{generic}
we know that finite convergence of Lasserre's LMI hierarchy
is ensured generically, yet we do not know a priori at which relaxation order it occurs.
To certify finite convergence, we can use the following condition.

\begin{proposition}[Certificate of finite convergence]\label{flat}
Let $y^*$ be the solution of LMI problem (\ref{hierarchy})
at a given relaxation order $r\geq r_X$.
If
\[
\mathrm{rank}\:M_{r-{r_X}}(y^*) = \mathrm{rank}\:M_r(y^*)
\]
then $p^*_r=p^*$.
\end{proposition}

If the moment matrix rank conditions of Proposition \ref{flat}
are satisfied, then we can use numerical linear algebra
to extract $\mathrm{rank}\:M_r(y^*)$ global optima for problem (\ref{poly}).
However, we do not describe the algorithm in this document.

\begin{proposition}[Rank-one moment matrix]\label{rankone}
The condition of Proposition \ref{flat}
is satisfied if
\[
\mathrm{rank}\:M_r(y^*) = 1.
\]
\end{proposition}

If the rank condition of Proposition \ref{rankone} is satisfied,
first order moments readily yield a global optimum:
$x^*=(y^*_{\alpha})_{|\alpha|=1}$.

\begin{example}
For the polynomial optimization problem
of Example \ref{ex:polyopt}, we obtain at the
second LMI relaxation a rank-one matrix $M_2(y^*)\geq 0$ and the
global optimum $x^*_1=y^*_{10}=\frac{1-\sqrt{5}}{2}$,
$x^*_2=y^*_{01}=\frac{1+\sqrt{5}}{2}$.
\end{example}

\section{Complexity estimates}

Consider a polynomial optimization problem 
\[
\begin{array}{lll}
p^* = & \min & p_0(x) \\
& \mathrm{s.t.} & p_k(x) \geq 0, \: k=1,\ldots,n_X
\end{array}
\]
as in (\ref{poly}), with $x \in {\mathbb R}^n$, and its hierarchy
of LMI relaxations (\ref{hierarchy}).

Let us denote by $M$ the number of variables, i.e. the size of vector $y$,
in the LMI relaxation of order $r$. It is equal to the number of monomials
of $n$ variables of degree $2r$, namely $M={n+2r \choose n}$.
If the number of variables $n$ is fixed (e.g. for
a given polynomial optimization problem) then $M$ grows in $O(r^n$),
that is polynomially in the relaxation order $r$.
If the relaxation order $r$ is fixed (say to the smallest possible
value, the first LMI relaxation in the hierarchy),
then $M$ grows in $O(n^r)$, that is polynomially in the number
of variables $n$.

In practice, given the current state-of-the-art in general-purpose
SDP solvers and personal computers, we can expect an LMI problem
to be solved in a matter of a few minutes provided the problem
is reasonably well-conditioned and the number of variables $M$
is less than $5000$, say.

\section{Convex hulls of semialgebraic sets}

Let us use the notations defined in Section \ref{sec:hierarchy},
and let $X$ be the compact basic semialgebraic set defined there
and satisfying Assumption \ref{compact}.
For $r \geq r_X$ consider the spectrahedral shadow
\[
X_d := \{(y_{\alpha})_{|\alpha|=1} \in {\mathbb R}^n \: :\:
y_0 = 1, \: M_r(y) \geq 0, \: M_{r-r_k}(p_k\:y) \geq 0, \: k=1,\ldots,n_X\}.
\]

\begin{proposition}[Convex outer approximations of semialgebraic sets]\label{outer}
$X_r$ is an outer approximation of $X$, i.e.
$X \subset X_r$. Moreover, $X_{r+1} \subset X_r$, and
$X_{\infty} = \mathrm{conv}\:X$.
\end{proposition}

The result is also true if $X$ is a compact algebraic set, defined
by finitely many polynomial equations. In particular, if $X$
is finite-dimensional, i.e. the union of a finite number of
points of ${\mathbb R}^n$, then $X_{\infty}$ is a polytope.

\begin{figure}
\begin{center}
\includegraphics[width=0.7\textwidth]{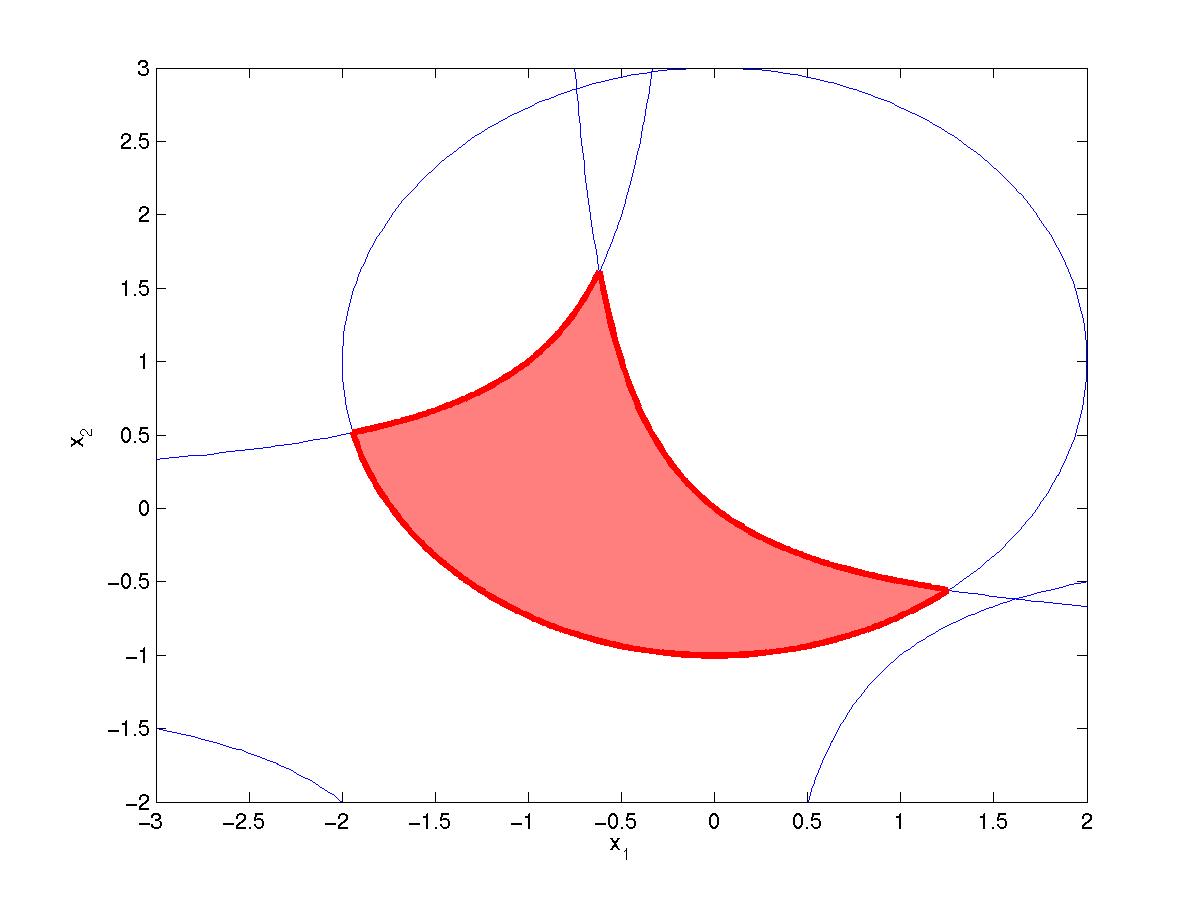}   
\caption{Nonconvex semialgebraic set $X$ (in red).\label{fig:polyopt1}}                              
\end{center}      
\end{figure}
\begin{figure}
\begin{center}
\includegraphics[width=0.7\textwidth]{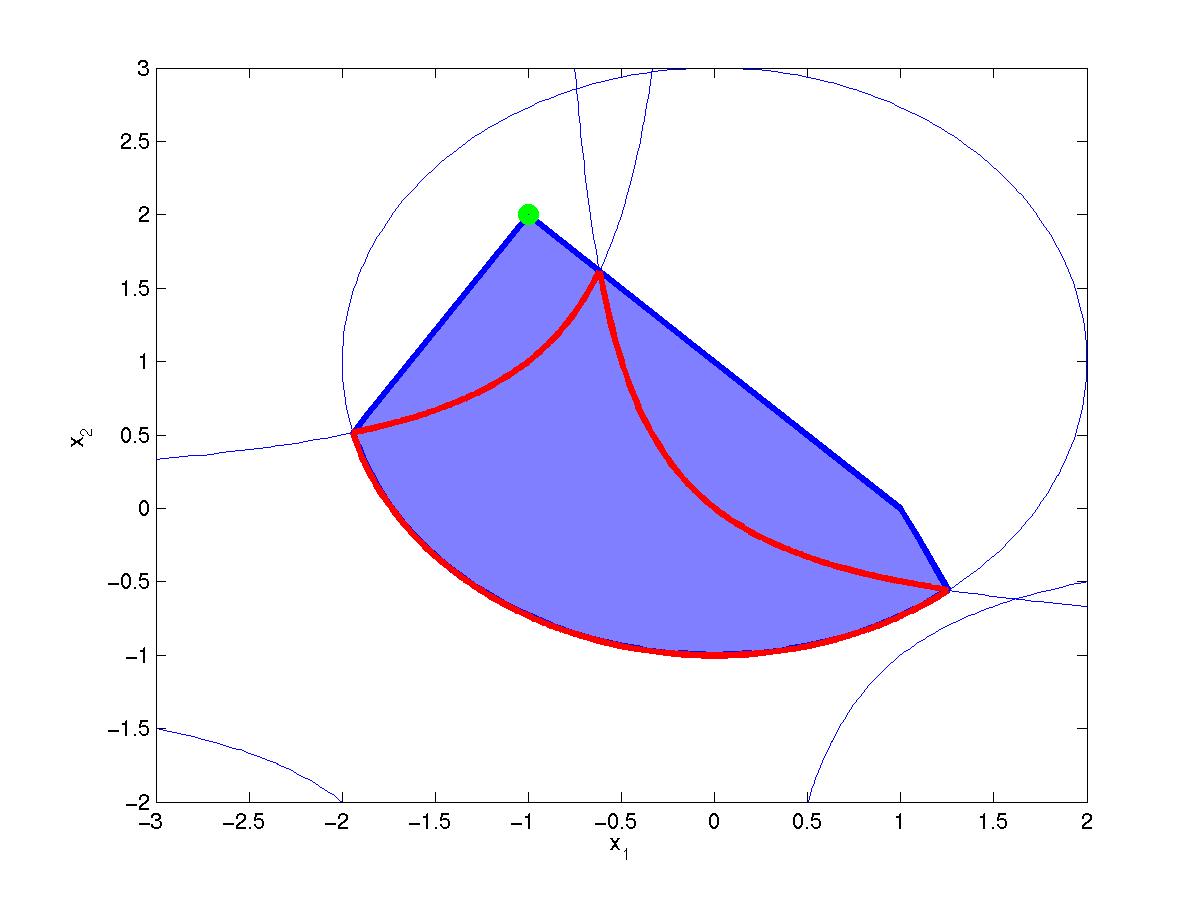}   
\caption{First spectrahedral shadow $X_1 \supset X$ (in blue)
with boundary of $X$ (in red) and suboptimal point (in green).
\label{fig:polyopt2}}                              
\end{center}      
\end{figure}
\begin{figure}
\begin{center}
\includegraphics[width=0.7\textwidth]{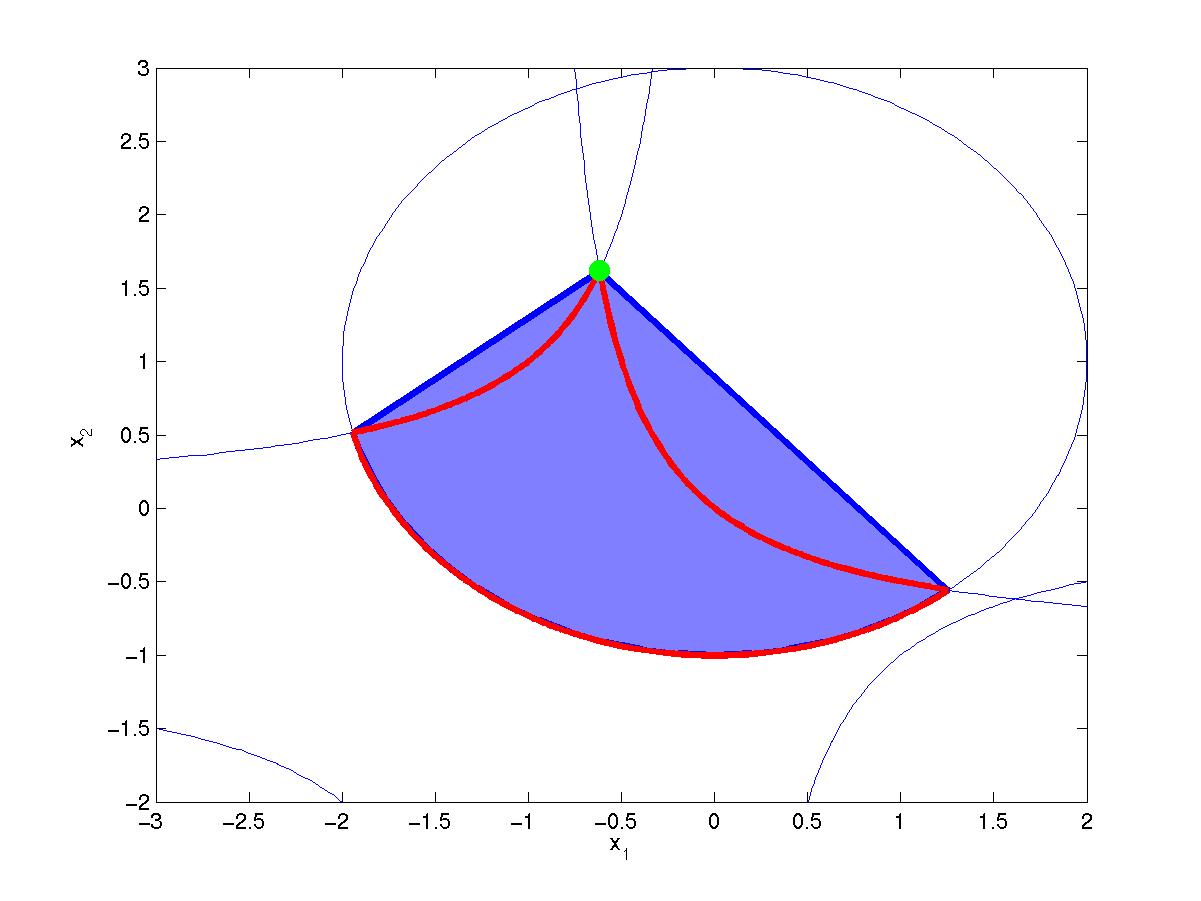}   
\caption{Second spectrahedral shadow $X_2 \supset X$ (in blue)
with boundary of $X$ (in red) and optimal point (in green).
\label{fig:polyopt3}}                           
\end{center}      
\end{figure}

\begin{example}\label{ex:polyopt2}
The polynomial optimization problem of Example \ref{ex:polyopt}
has a compact basic semialgebraic feasible set
\[
X=\{x \in {\mathbb R}^2 \: :\: 3-2x_2+x^2_1-x^2_2 \geq 0, \:
-x_1-x_2-x_1x_2 \geq 0, \: 1+x_1x_2 \geq 0\}
\]
represented in red on Figure \ref{fig:polyopt1}.
The first spectrahedral shadow $X_1 \supset X$, corresponding to the projection
on the plane of first-order moments of the 5-dimensional spectrahedron
of the first LMI relaxation, is represented in blue on Figure \ref{fig:polyopt2}.
Also represented in green is the point corresponding to the minimization
of $-y_{01}$, yielding the lower bound $p^*_1$.
The second spectrahedral shadow $X_2 \supset X$, corresponding to the projection
of the plane of first-order moments of the 14-dimensional spectrahedron
of the second LMI relaxation, is represented in blue on Figure \ref{fig:polyopt3}.
Also represented in green is the point corresponding to the minimization
of $-y_{01}$, yielding the lower bound $p^*_2$.
Apparently, $X_2=\mathrm{conv}\:X$, so minimizing $-x_2$ on $X$
or $-y_{01}$ on $X_2$ makes no difference.
\end{example}

\section{Software interfaces}

A {\tt Matlab} interface called {\tt GloptiPoly} has been
designed to construct Lasserre's LMI relaxations in a format
understandable by the SDP solver {\tt SeDuMi}, but also
any other SDP solver interfaced via {\tt YALMIP}.
It can be used to construct an LMI relaxation (\ref{hierarchy})
of given order corresponding to a polynomial optimization problem (\ref{poly})
with given polynomial data entered symbolically.
More generally, it can be used to model generalized problems
of moments (\ref{pmom}). A numerical algorithm is implemented
in {\tt GloptiPoly} to detect global optimality of an LMI
relaxation, using the rank tests of Propositions \ref{flat}
and \ref{rankone}. The algorithm also extracts
numerically the global optima from a singular value
decomposition of the moment matrix.
Another {\tt Matlab} interface called {\tt SOSTOOLS} was
developed independently and concurrently. It focuses on
the dual polynomial sum-of-squares decompositions mentioned
at the end of Section \ref{sec:momlp}, but not described
in this document. Note however that there is no global
optimality detection and global optima extraction algorithm
in {\tt SOSTOOLS}. Specialized {\tt moment} and {\tt sos}
modules are available in the interface {\tt YALMIP}
that implement some of the algorithms of {\tt GloptiPoly}
and {\tt SOSTOOLS}. For sparse polynomial optimization
problems (with polynomial data featuring a few number
of nonzero monomials), a specialized interface called
{\tt SparsePOP} is available. It generates reduced-size
LMI relaxations by exploiting the problem structure.

Note that these interfaces only generate the LMI relaxations
in a format understandable by general-purpose SDP solvers.
There is currently no working implementation of a
dedicated SDP solver for problems coming from
polynomial optimization.

\section{Back to the motivating example}\label{eig2}

Let us address the eigenvalue assignment problem
of Section \ref{eig}. We formulate it as a nonconvex
polynomial optimization problem
\[
\begin{array}{ll}
\min & p_0(x) \\
\mathrm{s.t.} & p_k(x) = 0, \:k=1,\ldots,n
\end{array}
\]
where the objective function is a positive definite convex
quadratic form
\[
p_0(x):=\sum_{i,j=1}^n (x_i-x_j)^2.
\]
This choice is motivated by physical reasons, and it
corresponds to the search of a solution $x$ with
entries $x_i$ as identical as possible.

First, let us generate the
system of polynomial equations $p_k(x)=0$, $k=1,\ldots,n$
with the following Maple script:
\begin{verbatim}
with(LinearAlgebra):with(PolynomialTools):
n:=3:B:=Matrix(n):for i from 1 to n-1 do
 B(i,i):=2: B(i,i+1):=-1: B(i+1,i):=-1:
end do: B(n,n):=(n+1)/n;
K:=Matrix(n,Vector(n,symbol=k),shape=diagonal):
q:=product(x-1/((2*j)^2-1),j=1..n):
p:=CoefficientList(collect(charpoly(MatrixInverse(B).F,x)-q,x),x);
\end{verbatim}
For $n=3$ this code generates the following polynomials
\[
\begin{array}{rcl}
p_1(x) & = & \frac{5}{6}x_1+\frac{4}{3}x_2+\frac{3}{2}x_3 -\frac{3}{7}\\
p_2(x) & = & \frac{2}{3}x_1x_2+x_1x_3+x_2x_3 -\frac{53}{1575}\\
p_3(x) & = & \frac{1}{2}x_1x_2x_3-\frac{1}{1575}.
\end{array}
\]
These polynomials are then converted into {\tt Matlab} format,
and we use the following {\tt GloptiPoly} code for inputing problem 
and solving the smallest possible LMI relaxation, i.e. $r=2$
in problem (\ref{hierarchy}):
\begin{verbatim}
mpol x 3
X = [5/6*x(1)+4/3*x(2)+3/2*x(3)-3/7
     2/3*x(1)*x(2)+x(1)*x(3)+x(2)*x(3)-53/1575
     1/2*x(1)*x(2)*x(3)-1/1575];
obj = 0;
for i = 1:length(x)
 for j = 1:length(x)
  obj = obj+(x(i)-x(j))^2;
 end
end
P = msdp(min(obj),X==0);
[stat,obj] = msol(P);
double(x)
\end{verbatim}
With this code and the SDP solver {\tt SeDuMi}, we obtain
the unique solution
\[
x \approx (9.3786\cdot10^{-2}, \; 8.6296\cdot10^{-2},\;
1.5690\cdot10^{-1})
\]
(to 5 significant digits) certified numerically by a rank one moment matrix,
see Proposition \ref{rankone},
after less than one second of CPU time on a standard desktop computer.

The cases $n=2,3,4,5$ are solved very
easily (in a few seconds) but the solution
(obtained with {\tt SeDuMi}) is not very accurate.
We have not investigated the possibility of refining
the solution with e.g. Newton's method. We have not
investigated either the possibility of certifying
rigorously the solution using e.g. {\tt VSDP} or
multiprecision arithmetic.

The case $n=6$ is solved in a few minutes, and the case $n=7$
is significantly harder: it takes a few hours to be solved.
Finally, solving the case $n=8$ takes approximately $15$ hours
on our computer.

\section{Notes and references}

Accessible introductions to measure theory and relevant notions
of functional analysis and probability theory are \cite{rn55,kf68,l69,rf10}.
Lasserre's hierarchy of LMI relaxations for polynomial optimization 
were originally proposed in \cite{l00,l01} with a proof of
convergence (Proposition \ref{lasserre}) relying on Putinar's
Positivstellensatz (Proposition \ref{putinar}) described
in \cite{p93}. The genericity result of Proposition \ref{generic}
is described in \cite{n12}. Examples \ref{ex:polyopt} and
\ref{ex:polyopt2} are from \cite{hl04}.
The rank condition of Proposition \ref{flat}
relies on flat extension results by Curto and Fialkow, see
e.g. \cite{cf05}. The algorithm implemented
in {\tt GloptiPoly} for extracting global optima was
described in \cite{hl05}, see also \cite{l09} and \cite{mom09}
for more comprensive descriptions. The use of dual
polynomial sum-of-squares was proposed in \cite{p00}, see
also \cite{bpt13}. Optimization over polynomial sum-of-squares
and more generally squared functional systems and the connection
with SDP was studied in \cite{n00}, mostly in the univariate case.
An excellent survey of this material (on both the primal problems
on moments and the dual problems on polynomial sum-of-squares)
is \cite{l09}, and the reader
is refered to \cite{mom09} for a more advanced treatment. 
The outer approximation result
of Proposition \ref{outer} follows from the convergence
proof of Lasserre's relaxation \cite{l01} and elementary duality
arguments.
Finally, the structured eigenvalue assignment problem of Sections \ref{eig}
and \ref{eig2} is comprehensively described in \cite{marx}.


\chapter{Infinite-dimensional polynomial optimization}\label{chap:ido}

We extend the approach of the previous chapter
to optimization over the infinite-dimensional sets
of solutions of ordinary differential equations
with polynomial vector fields.

\section{Occupation measures}\label{occmea}

Let us consider the nonlinear ordinary differential equation (ODE)
\begin{equation}\label{ode}
\dot{x}(t) = f(t,x(t))
\end{equation}
for $t \in [0,T]$ with a given terminal time $T>0$,
where $x : [0,T] \to {\mathbb R}^n$ is a time-dependent $n$-dimensional state vector,
and vector field $f : [0,T]\times {\mathbb R}^n \to {\mathbb R}^n$ is a smooth map.
Given a set $X \subset {\mathbb R}^n$, we assume that dynamics $f$ and terminal
time $T$ are such that there is a solution to the Cauchy problem for
ODE (\ref{ode}). Since vector field $f$ is smooth, this solution is unique
for any given initial condition $x(0)=x_0 \in X$. Any such solution, or trajectory, $x(t)$
is an absolutely continuous function of time with values in $X$, and to emphasize the
dependence of the solution on the initial condition we write $x(t\:|\:x_0)$.

Now think of initial condition $x_0$ as a random variable in $X$, or more
abstractly as a probability measure $\xi_0 \in {\mathscr P}(X)$,
that is a map from the Borel $\sigma$-algebra 
${\mathscr B}(X)$ of subsets of $X$ to the interval $[0,1] \subset {\mathbb R}$
such that $\xi_0(X)=1$. For example, the expected value of $x_0$ is the vector
$E[x_0] = \int_X x\:\xi_0(dx)$ of first order moments of $\xi_0$.

Now solve ODE (\ref{ode}) for a trajectory, given this
random initial condition. At each time $t$, the state can also be interpreted
as a random variable, i.e. a probability measure that we denote by $\xi \in {\mathscr P}(X)$.
We say that the measure is transported by the flow of the ODE.
We also use the notation $\xi(dx\:|\:t)$ if we want to emphasize the fact that $\xi$
is a conditional probability measure, or stochastic kernel,
i.e. a probability measure acting on subsets of
${\mathscr B}(X)$ for each given, or frozen value of $t$.

This one-dimensional family, or path of measures, satisfies a partial differential equation
(PDE) which turns out to be linear in the space of probability measures.
This PDE is usually called Liouville's equation.
Conversely, the nonlinear ODE follows by applying Cauchy's method of characteristics
to the linear transport PDE.

Let us now derive the Liouville equation explicitly.

\begin{definition}[Indicator function]
The indicator function of a set A is the function $x\mapsto I_A(x)$
such that $I_A(x)=1$ when $x\in A$ and $I_A(x)=0$ when $x\notin A$.
\end{definition}

\begin{definition}[Occupation measure]
Given an initial condition $x_0$, the occupation measure
of a trajectory $x(t\:|\:x_0)$ is defined by
\[
\mu(A\times B\:|\:x_0):=\int_A I_B(x(t\:|\:x_0))dt
\]
for all $A \in {\mathscr B}([0,T])$ and $B \in {\mathscr B}(X)$.
\end{definition}

A geometric interpretation is that $\mu$ measures the time spent by the graph
of the trajectory $(t,x(t\:|\:x_0))$ in a given subset $A\times B$ of
$[0,T]\times X$. An analytic interpretation is that integration w.r.t $\mu$ 
is equivalent to time-integration along a system trajectory, i.e.
\[
\int_0^T v(t,x(t\:|\:x_0))dt = \int_0^T \int_X v(t,x)\mu(dt,dx\:|\:x_0)
\]
for every test function $v \in {\mathscr C}([0,T]\times X)$.

\begin{example}[Occupation measure for a scalar linear system]\label{ex:occlin}
Consider the one-dimensional ODE $\dot{x}(t)=-x(t)$ with initial
condition $x(0)=x_0 \geq 0 $, whose solution is $x(t)=x_0e^{-t}$.
Given $a \geq 0$, the occupation measure of the trajectory is such that
\[
\begin{array}{rcll}
\mu([0,1]\times[0,a] \:|\:x_0)
& = & 1 & \mathrm{if}\: x_0 \leq a \\
& = & 1-\log\frac{x_0}{a} & \mathrm{if}\: a \leq x_0 \leq ae \\
& = & 0 & \mathrm{if}\: x_0 > ae \\
\end{array}
\]
where $e\approx 2.71828$ is Euler's number.
\end{example}

Now define the linear operator ${\mathcal L} : {\mathscr C}^1([0,T]\times X)
\to {\mathscr C}([0,T]\times X)$ by
\[
v \mapsto {\mathcal L}v := \frac{\partial v}{\partial t} + 
\sum_{i=1}^n \frac{\partial v}{\partial x_i} f_i =
\frac{\partial v}{\partial t} + (\mathrm{grad}\:v)' f
\]
and its adjoint operator ${\mathcal L}' : {\mathscr C}([0,T]\times X)' \to
{\mathscr C}^1([0,T]\times X)'$ by the relation
\[
\langle v,{\mathcal L}'\mu\rangle := \langle {\mathcal L}v,\mu\rangle
= \int_0^T \int_X {\mathcal L}v(t,x,u)\mu(dt,dx)
\]
for all $\mu \in {\mathscr M}([0,T]\times X)={\mathscr C}([0,T]\times X)'$
and $v \in {\mathscr C}^1([0,T]\times X)$. This operator can also be
expressed as
\[
\mu \mapsto {\mathcal L}'\mu = -\frac{\partial \mu}{\partial t}
-\sum_{i=1}^n \frac{\partial (f_i\mu)}{\partial x_i} =
-\frac{\partial \mu}{\partial t} - \mathrm{div}\:f\mu
\]
where the derivatives of measures are understood in the sense of distributions
(i.e. via their action on smooth test functions), and the change of sign comes
from the integation by parts formula.

Given a test function $v \in {\mathscr C}^1([0,T]\times X)$ it follows from
the above definition of the occupation measure that
\begin{equation}\label{intv}
\begin{array}{rcl}
v(T,x(T\:|\:x_0)) & = & v(0,x_0) + \int_0^T \dot{v}(t,x(t\:|\:x_0))dt \\
& = & v(0,x_0) + \int_0^T {\mathcal L}v(t,x(t\:|\:x_0))dt \\
& = & v(0,x_0) + \int_0^T \int_X {\mathcal L}v(t,x)\mu(dt,dx\:|\:x_0).
\end{array}
\end{equation}

\begin{definition}[Initial measure]\label{initial}
The initial measure $\xi_0 \in {\mathscr P}(X)$ is a probability measure
that rules the distribution in space of the initial condition $x_0$.
\end{definition}

\begin{definition}[Average occupation measure]
Given an initial measure $\xi_0$, the average occupation measure
of the flow of trajectories is defined by
\[
\mu(A\times B):=\int_X \mu(A\times B\:|\:x_0)\xi_0(dx_0)
\]
for all $A \in {\mathscr B}([0,T])$ and $B \in {\mathscr B}(X)$.
\end{definition}

\begin{example}
Returning to the scalar linear ODE of Example \ref{ex:occlin},
with initial conditions uniformly distributed on $[0,1]$, i.e.
$\xi_0(dx)=I_{[0,1]}(x)dx$, the average occupation measure
is such that
\[
\mu([0,1]\times[0,a])=\int_0^1 \mu([0,1]\times[0,a]\:|\:x_0)dx_0 = \int_0^a dx_0
+ \int_a^{ae} \left(1-\log\frac{x_0}{a}\right)dx_0 = a(e-1)
\]
for any given $a\geq 0$.
\end{example}

\begin{definition}[Terminal measure]\label{terminal}
The terminal measure $\xi_T \in {\mathscr P}(X)$ is a probability measure
that rules the distribution in space of the terminal condition $x(T)$.
It is defined by
\[
\xi_T(B):=\int_X I_B(x(T\:|\:x_0))\xi_0(dx_0)
\]
for all $B \in {\mathscr B}(X)$.
\end{definition}

It follows by integrating equation (\ref{intv}) with respect to $\xi_0$
that
\[
\int_X v(T,x)\xi_T(dx) = \int_X v(0,x)\xi_0(dx) +
\int_0^T \int_X {\mathcal L}v(t,x)\mu(dt,dx)
\]
or more concisely
\begin{equation}\label{liouv}
\langle v(T,.),\xi_T \rangle = \langle v(0,.),\xi_0 \rangle + \langle {\mathcal L}v,\mu \rangle
\end{equation}
which is a linear equation linking the initial measure $\xi_0$,
the terminal measure $\xi_T$ and the occupation measure $\mu$,
for all $v \in {\mathscr C}^1([0,T]\times X)$.

Letting
\[
\mu_0(dt,dx) := \delta_0(dt)\:\xi_0(dx), \quad
\mu_T(dt,dx) := \delta_T(dt)\:\xi_T(dx),
\]
we can write $\langle v(0,.),\xi_0 \rangle =
\langle v, \mu_0 \rangle$ and $\langle v(T,.),\xi_T \rangle =
\langle v, \mu_T \rangle$. Then, equation (\ref{liouv}) can
be rewritten equivalently using the adjoint linear operator as
\[
\langle v, {\mathcal L}'\mu  \rangle = \langle v, \mu_T \rangle
- \langle v, \mu_0 \rangle
\]
and since this equation is required to hold for all test functions
$v \in {\mathscr C}^1([0,T]\times X)$, we obtain
a linear PDE on measures ${\mathcal L}'\mu=\mu_T-\mu_0$ that we write
\begin{equation}\label{liouville}
\frac{\partial \mu}{\partial t}+\mathrm{div}\:f\mu = \mu_0 - \mu_T
\end{equation}
where the derivatives should be understood in the sense of distributions.
This equation is classical in fluid mechanics and statistical physics, and
it is called the equation of conservation of mass, or the continuity equation,
or the advection equation, or {\bf Liouville's equation}.

Note that we can disintegrate the average occupation measure as follows
\[
\mu(dt,dx) = dt\:\xi(dx\:|\:t)
\]
where $\xi(.\:|\:t) \in {\mathscr P}(X)$
is the conditional of $\mu$ w.r.t. $t$,
and $dt$ is the marginal of $\mu$ w.r.t. $t$, here the Lebesgue measure.
Liouville's equation (\ref{liouville}) can be also written as
a linear PDE satisfied by probability measure $\xi$, namely
\begin{equation}\label{liouvillexi}
\frac{\partial \xi}{\partial t}+\mathrm{div}\:f\xi = 0
\end{equation}
with a given initial measure $\xi(.\:|\:t=0)=\xi_0$.

\begin{proposition}\label{continuity}
Given $\xi_0 \in {\mathscr P}(X)$, there is a unique solution $\xi(.|t) \in {\mathscr P}(X)$
solving equation (\ref{liouvillexi}). Letting $\xi_T:=\xi(.|t=T) \in {\mathscr P}(X)$,
there is a unique solution $\mu \in {\mathscr M}_+([0,T]\times X)$ solving equation (\ref{liouville}).
\end{proposition}

Note that in particular if $\xi_0 = \delta_{x_0}$, then $\xi(.|t)=\delta_{x(t\:|\:x_0)}$
is the Dirac measure supported on the trajectory $x(t\:|\:x_0)$ starting from $x_0$.
The geometric picture behind Liouville's equation (\ref{liouville})
is that it encodes a superposition of all classical solutions
solving the Cauchy problem (\ref{ode}). 
The main advantage of the Liouville PDE is that it is a linear equation
(in the infinite-dimensional space of measures), whereas the original
Cauchy ODE is nonlinear (in the infinite-dimensional space of
absolutely continuous trajectories).

\section{Measure LP}

Consider now the following dynamic optimization problem with polynomial differential constraints
\begin{equation}\label{optode}
\begin{array}{lll}
p^* \:= & \inf & \int_0^T l(t,x(t))dt \\
& \mathrm{s.t.} & \dot{x}(t) = f(t,x(t)), \:\: x(t) \in X, \:\:t \in [0,T] \\
& & x(0) \in {X}_0, \:\: x(T) \in {X}_T
\end{array}
\end{equation}
with given polynomial dynamics $f \in {\mathbb R}[t,x]$ and Lagrangian
$l \in {\mathbb R}[t,x]$, and state trajectory $x(t)$
constrained in a compact basic semialgebraic set
\[
{X} = \{x \in {\mathbb R}^n \: : \: p_k(t,x) \geq 0, \:
k=1,\ldots,n_X\}
\]
for given polynomials $p_k \in {\mathbb R}[t,x]$.
Finally, initial and terminal states are constrained in compact
basic semialgebraic sets
\[
{X}_0 = \{x \in {\mathbb R}^n \: : \: p_{0k}(x) \geq 0, \:
k=1,\ldots,n_0\} \subset {X}
\]
and
\[
{X}_T = \{x \in {\mathbb R}^n \: : \: p_{Tk}(x) \geq 0, \:
k=1,\ldots,n_T\} \subset {X}
\]
for given polynomials $p_{0k}, p_{Tk} \in {\mathbb R}[x]$.
In problem (\ref{optode}) the infimum is w.r.t. a trajectory $x(t)$
starting in $X_0$, ending in $X_T$, and staying in $X$.

Using the framework described in Section \ref{occmea},
we encode the state trajectory $x(t)$ in an occupation measure $\mu$
and we come up with an infinite-dimensional LP problem
\begin{equation}\label{pdolp}
\begin{array}{lll}
p^* \:= & \inf & \int l\mu \\
& \mathrm{s.t.} & \int \left( \frac{\partial v}{\partial t} +
(\mathrm{grad}\:v)'f \right)\:\mu  =
\int v \mu_T  - \int v \mu_0 \\
\end{array}
\end{equation}
for all smooth test functions $v \in {\mathscr C}^1([0,T]\times X)$
and where the infimum is w.r.t. occupation measure $\mu \in {\mathscr M}_+([0,T]\times X)$,
initial measure $\mu_0 \in {\mathscr P}(\{0\}\times X_0)$,
terminal measure $\mu_T \in {\mathscr P}(\{T\}\times X_T)$, and terminal time $T$.
Note that $\mu_0$ resp. $\mu_T$ and $T$ can be free, or given.
More abstractly, problem (\ref{pdolp})
can be written as a measure LP
\begin{equation}\label{lpmu0}
\begin{array}{lll}
p^* \:= & \inf & \langle l,\mu \rangle \\
& \mathrm{s.t.} & \frac{\partial \mu}{\partial t} + \mathrm{div}\:f\mu = \mu_0 - \mu_T\\
\end{array}
\end{equation}
where the linear constraint is Liouville's equation, and the minimum
is w.r.t. measures $(\mu,\:\mu_0,\:\mu_T) \in {\mathscr M}_+([0,T]\times X)\times
{\mathscr M}_+(\{0\}\times X_0)\times{\mathscr M}_+(\{T\}\times X_T)$.
If the three measures $\mu$, $\mu_0$ and $\mu_T$ are unknown, then
an additional linear constraint like $\mu_0(\{0\}\times X_0)=1$
or $\mu_T(\{T\}\times X_T)=1$ must be enforced
to rule out the trivial zero solution.

\begin{remark}[Autonomous case]\label{autonomous}
If the terminal time $T$ is free and the Lagrangian $l$ and the dynamics
$f$ do not depend explicitly on time $t$, then it can be shown
without loss of generality that in problem (\ref{lpmu0})
the measures do not depend explicitly on time either, and
the terminal time is equal to the mass of the occupation measure, i.e. $T = \mu(X)$.
The measure LP becomes
\begin{equation}\label{lpmu0a}
\begin{array}{lll}
p^* \:= & \inf & \langle l,\mu \rangle \\
& \mathrm{s.t.} & \mathrm{div}\:f\mu = \mu_0 - \mu_T\\
\end{array}
\end{equation}
where the minimum is taken w.r.t.
$(\mu,\:\mu_0,\:\mu_T) \in {\mathscr M}_+(X)\times
{\mathscr M}_+(X_0)\times{\mathscr M}_+(X_T)$.
\end{remark}

\begin{example}\label{ex:occtraj}
Consider again the scalar linear ODE of Example \ref{ex:occlin}
\[
\dot{x}=-x
\]
with initial measure $\mu_0(dt,dx):=\delta_{0}(dt)\:\xi_0(dx)$
with state distribution $\xi_0 \in {\mathscr P}(X_0)$ supported on 
\[
X_0:=\{x \in {\mathbb R} \: :\: p_0(x):=\frac{1}{4}-\left(x-\frac{3}{2}\right)^2\geq 0\},
\]
with terminal measure $\mu_T(dt,dx):=\delta_{T}(dt)\:\xi_T(dx)$
with state distribution $\xi_T \in {\mathscr P}(X_T)$ supported on
\[
X_T:=\{x \in {\mathbb R} \: :\: p_T(x):=\frac{1}{4}-x^2 \geq 0\},
\]
and with average occupation measure $\mu(dt,dx) := dt\:\xi(dx\:|\:t)$
with state conditional $\xi(dx\:|\:t) \in {\mathscr P}(X)$ supported for each $t \in [0,T]$ on
\[
X:=\{x \in {\mathbb R} \: :\: p(x):=4-x^2 \geq 0\}.
\]
We want to find trajectories minimizing the state energy $\int_0^T x^2(t)dt$.

The linear measure problem (\ref{lpmu0}) reads
\[
\begin{array}{lll}
p^* \:= & \inf & \langle x^2,\mu \rangle \\
& \mathrm{s.t.} & \frac{\partial \mu}{\partial t} - \mathrm{div}\:x\mu = \mu_0 - \mu_T\\
\end{array}
\]
where the minimum is w.r.t. terminal time $T$ and
nonnegative measures
$\mu$, $\mu_0$ and $\mu_T$ supported respectively
on $[0,T]\times X$, $\{0\}\times X_0$ and $\{T\}\times X_T$,
and we have to enforce the additional normalization constraint
$\mu_0(\{0\}\times X_0) = 1$.

This problem can be solved analytically, with optimal trajectory
$x(t)=e^{-t}$ leaving $X_0$ at $x(0)=1$ and reaching $X_T$ at $x(T)=\frac{1}{2}$
for $T=\log 2\approx 0.6931$. So the optimal measures solving the above
LP are
\[
\mu(dt,dx) = dt\:\delta_{e^{-t}}(dx), \quad
\mu_0(dt,dx) = \delta_{0}(dt) \: \delta_{1}(dx), \quad
\mu_T(dt,dx) = \delta_{\log 2}(dt) \: \delta_{\frac{1}{2}}(dx)
\]
and $p^*=\int_0^{\log 2}e^{-2t}dt = \frac{3}{8}$.

Alternatively, following Remark \ref{autonomous}, since the trajectory
optimization problem is autonomous,
we can also formulate the measure LP problem 
\[
\begin{array}{lll}
p^* \:= & \inf & \langle x^2,\mu \rangle \\
& \mathrm{s.t.} & - \mathrm{div}\:x\mu = \mu_0 - \mu_T\\
\end{array}
\]
w.r.t.
nonnegative measures $\mu$, $\mu_0$ and $\mu_T$ supported respectively
on $X$, $X_0$ and $X_T$, and the optimal solution of
the problem is now
\[
\mu(dx) = \int_0^T \delta_{e^{-t}}(dx)dt, \quad
\mu_0(dx) = \delta_{1}(dx), \quad
\mu_T(dx) = \delta_{\frac{1}{2}}(dx).
\]
\end{example}

\section{Moment LP and LMI relaxations}\label{inflplmi}

Let us write problem (\ref{lpmu0}) as
as special instance of a more general measure LP
\[
\begin{array}{rcll}
p^* & = & \inf & \langle c,\nu \rangle \\
&& \mathrm{s.t.} & {\mathcal A}\nu = \beta \\
&&& \nu \in {\mathscr M}^n_+
\end{array}
\]
where the decision variable is 
an $n$-dimensional vector of nonnegative measures $\nu$.
Linear operator ${\mathcal A} : {\mathscr M}^n \to {\mathscr M}^m$ takes
an $n$-dimensional vector of measures and returns an $m$-dimensional vector
of measures. The right hand side $\beta \in {\mathscr M}^m$ is a given
$m$-dimensional vector of measures. The objective function is the
duality pairing between a given $n$-dimensional vector of continuous functions
$c \in {\mathscr C}^n$ and $\nu$,
i.e. $\langle c,\nu \rangle = \sum_{i=1}^n \langle c_i,\nu_i \rangle =
\sum_{i=1}^n \int c_i \nu_i $. If we suppose that all the functions
are polynomials, i.e. $a_{ij}(x) \in {\mathbb R}[x]$,
$c_i(x) \in {\mathbb R}[x]$, $i=1,\ldots,n$, $j=1,\ldots,m$,
then each measure $\nu_i$ can be manipulated via the sequence
$y_i:=(y_{i\alpha})_{\alpha \in {\mathbb N}}$ of its moments.
The measure LP becomes a moment LP
\begin{equation}\label{pmominf}
\begin{array}{rcll}
p^* & = & \inf & \sum_{i=1}^n \sum_{\alpha} c_{i\alpha}y_{i\alpha} \\
&& \mathrm{s.t.} & \sum_{i=1}^n \sum_{\alpha} a_{ij\alpha}y_{i\alpha} = b_j, \:\: j=1,\ldots,m \\
&&& \text{$y_i$ has a representing measure $\nu_i \in {\mathscr M}_+(X_i)$, $i=1,\ldots,n$.}
\end{array}
\end{equation}
As in Section \ref{sec:hierarchy}, we use the explicit LMI conditions of Section \ref{sec:repres}
to model the constraints that a sequence has a representing measure.
If each semialgebraic set
\[
X_i := \{x \in {\mathbb R}^n \: :\: p_{ik}(x) \geq 0, \: k=1,\ldots,n_i\}
\]
satisfies Assumption \ref{compact} for $i=1,\ldots,n$, problem (\ref{pmominf}) becomes 
\[
\begin{array}{rcll}
p^* & = & \inf & \sum_{i=1}^n \sum_{\alpha} c_{i\alpha}y_{i\alpha} \\
&& \mathrm{s.t.} & \sum_{i=1}^n \sum_{\alpha} a_{ij\alpha}y_{i\alpha} = b_j, \:\: j=1,\ldots,m \\
&&& M(y_i) \geq 0, \: M(p_{ik}\:y_i) \geq 0, \: i=1,\ldots,n, \: k=1,\ldots,n_i.
\end{array}
\]
Then, in order to solve problem (\ref{lpmu0}), we can build a hierarchy of
finite-dimensional LMI relaxations.
This generates a monotonically nondecreasing 
sequence of lower bounds asymptotically converging to $p^*$.
Details are omitted. 

\begin{example}\label{ex:occtraj2}
At the end of Example \ref{ex:occtraj} we came up with the autonomous measure LP
\[
\begin{array}{lll}
p^* \:= & \inf & \langle x^2,\mu \rangle \\
& \mathrm{s.t.} & - \mathrm{div}\:x\mu = \mu_0 - \mu_T\\
\end{array}
\]
in the decision variables $(\mu,\:\mu_0,\:\mu_T) \in {\mathscr M}_+(X)\times{\mathscr M}_+(X_0)
\times{\mathscr M}_+(X_T)$ with normalization constraint $\mu_0(X_0)=1$.
The corresponding moment LP problem (\ref{pmominf})  reads
\[
\begin{array}{rcll}
p^* & = & \inf & \int x^2 \mu(dx) \\
&& \mathrm{s.t.} & \int \mu_0(dx) = 1 \\
&&& -\alpha \int x^{\alpha} \mu(dx) = \int x^{\alpha} \mu_T(dx) - \int x^{\alpha} \mu_0(dx),
\quad \alpha = 0,1,2,\ldots
\end{array}
\]
or equivalently
\[
\begin{array}{rcll}
p^* & = & \inf & y_2 \\
&& \mathrm{s.t.} & {y_0}_0 = 1 \\
&&& -\alpha y_{\alpha} = {y_T}_{\alpha} - {y_0}_{\alpha}, \:\: \alpha=0,1,2\ldots \\
&&& \text{$y$ has a representing measure $\mu \in {\mathscr M}_+(X)$} \\
&&& \text{$y_0$ has a representing measure $\mu_0 \in {\mathscr M}_+(X_0)$} \\
&&& \text{$y_T$ has a representing measure $\mu_T \in {\mathscr M}_+(X_T)$} \\
\end{array}
\]
and the corresponding LMI relaxation of order $r$ is given by
\[
\begin{array}{rcll}
p^*_r & = & \inf & y_2 \\
&& \mathrm{s.t.} & {y_0}_0 = 1 \\
&&& -\alpha y_{\alpha} = {y_T}_{\alpha} - {y_0}_{\alpha}, \:\: \alpha=0,1,\ldots,2r \\
&&& M_r(y) \geq 0, \: M_{r-1}(p\:y) \geq 0\\
&&& M_r(y_0) \geq 0, \: M_{r-1}(p_0\:y_0) \geq 0\\
&&& M_r(y_T) \geq 0, \: M_{r-1}(p_T\:y_T) \geq 0
\end{array}.
\]
From the analytic solution described in Example \ref{ex:occtraj} we
can compute the entries of the moment vector $y$ of measure $\mu$,
namely $y_0=\log 2$ and
\[
y_{\alpha} = \int x^{\alpha}\mu(dx) = \int_0^{\log 2} e^{-\alpha t}dt = \frac{1-2^{-\alpha}}{\alpha},
\quad \alpha=1,2,\ldots
\]
\end{example}

\section{Optimal trajectory recovery}\label{sec:rectraj}

Once an LMI relaxation of given order is solved,
we expect vector $y$ to contain approximate moments
of the optimal occupation measure corresponding to
the optimal trajectory (if it is unique),
or at least a superposition (convex combination)
of optimal trajectories. In some cases, we can
recover approximately the trajectory from the
knowledge of its moments. The dual LMI relaxations
can be useful for this purpose. However, we do not elaborate
further on this point in this document.

\begin{example}
Solving the LMI relaxations of Example \ref{ex:occtraj2},
we observe that the moment matrices of the initial
and terminal measures both have rank one (to numerical
roundoff errors), with respective moment vectors
\[
{y_0}_{\alpha} = \int x^{\alpha} \mu_0(dx) = 1, \quad
{y_T}_{\alpha} = \int x^{\alpha} \mu_T(dx) = 2^{-\alpha}, \quad \alpha=0,1,2,\ldots
\]
From this it follows that $\mu_0=\delta_{1}$,
$\mu_T=\delta_{\frac{1}{2}}$ and the unique optimal trajectory
starts from $x(0)=1$ and reaches $x(T)=\frac{1}{2}$.
\end{example}

\section{Extension to piecewise polynomial dynamics}\label{sec:ppd}

The framework can be extended readily to differential equations
with terminal cost and piecewise polynomial dynamics
\begin{equation}\label{ppd}
\begin{array}{lll}
p^* \:= & \inf & f_0(x(T)) + \int_0^T l(t,x(t))dt \\
& \mathrm{s.t.} & \dot{x}(t) = f_j(t,x(t)), \:\: x(t) \in {X}_j, \:\: j=1,\ldots,N,
\:\: t \in [0,T] \\
& & x(0) \in {X}_0, \:\: x(T) \in {X}_T
\end{array}
\end{equation}
with given polynomial dynamics $f_j \in {\mathbb R}[t,x]$, Lagrangian
$l \in {\mathbb R}[t,x]$, terminal cost $f_0 \in {\mathbb R}[x]$ and state trajectory $x(t)$
constrained in compact basic semialgebraic sets $X_j$.
We assume that the state-space partitioning sets, or cells ${X}_j$, are disjoint,
i.e. all their respective intersections have zero Lebesgue measure in ${\mathbb R}^n$,
and they all belong to a given compact semialgebraic set $X$, e.g. a Euclidean ball of
large radius. Initial and terminal states are constrained in given compact basic
semialgebraic sets $X_0$ and $X_T$.

We can then extend the measure LP framework
to several measures $\mu_j$, one supported on each cell $X_j$, so
that the global occupation measure is
\[
\mu = \sum_{j=1}^N \mu_j.
\]
The measure LP reads
\[
\begin{array}{lll}
p^* \:= & \inf & \langle f_0,\mu_T \rangle + \sum_{j=1}^N \langle l,\mu_j \rangle \\
& \mathrm{s.t.} & \sum_{j=1}^N \left(\frac{\partial \mu_j}{\partial t} + 
\mathrm{div}\:f_j \mu_j \right) + \mu_T = \mu_0.
\end{array}
\]
It can solved numerically with a hiearchy of LMI relaxations
as in Section \ref{inflplmi}.

\section{Back to the motivating example}\label{valid2}

To address our motivating problem of Section \ref{valid}
we formulate an optimization problem (\ref{ppd}) with systems dynamics defined
as locally affine functions in three cells $X_j$, $j=1,2,3$
corresponding respectively to the linear regime of the torque saturation
\[
X_1 = \{x \in {\mathbb R}^2 \: :\: |F'x| \leq L\}, \quad f_1(x) = \left[\begin{array}{c}
x_1 \\ - F'x \end{array}\right]
\]
the upper saturation regime
\[
X_2 = \{x \in {\mathbb R}^2 \: :\: F'x \geq L\}, \quad f_2(x) = \left[\begin{array}{c}
x_1 \\ - L \end{array}\right]
\]
and the lower saturation regime
\[
X_3 = \{x \in {\mathbb R}^2 \: :\: F'x \leq -L\}, \quad f_3(x) = \left[\begin{array}{c}
x_1 \\ L \end{array}\right].
\]
The objective function has no integral term and a concave quadratic
terminal term $f_0(x) = -x(T)^Tx(T)$ which we would like to minimize,
so as to find trajectories with terminal states of largest norm.
If we can certify that for every initial state $x(0)$ chosen in $X_0$
the final state $x(T)$ belongs to set included in the
deadzone region, we have validated our controlled system.
In the measure LP problem of Section \ref{sec:ppd},
the 3 measures $\mu$, $\mu_0$ and $\mu_T$ are unknown
so we have to insert a normalization constraint to rule out the
trivial zero solutions:
\[
\begin{array}{lll}
p^* \:= & \inf & \langle f_0,\mu_T \rangle + \sum_{j=1}^N \langle l,\mu_j \rangle \\
& \mathrm{s.t.} & \sum_{j=1}^N \left(\frac{\partial \mu_j}{\partial t} + 
\mathrm{div}\:f_j \mu_j \right) = \mu_0 - \mu_T\\
&& \mu_0(\{0\}\times X_0) = 1.
\end{array}
\]

The resulting GloptiPoly script,
implementing some elementary scaling
strategies to improve numerical behavior of the SDP solver,
is as follows:

\begin{verbatim}
I = 27500; % inertia
kp = 2475; kd = 19800; % controller gains
L = 380; % input saturation level 
dz1 = 0.2*pi/180; dz2 = 0.05*pi/180; % deadzone levels
thetamax = 50; omegamax = 5; % bounds on initial conditions
epsilon = sqrt(1e-5); % bound on norm of terminal condition
T = 50; % final time

r = input('order of relaxation ='); r = 2*r;

% measures
mpol('x1',2); m1 = meas(x1); % linear regime
mpol('x2',2); m2 = meas(x2); % upper sat
mpol('x3',2); m3 = meas(x3); % lower sat
mpol('x0',2); m0 = meas(x0); % initial
mpol('xT',2); mT = meas(xT); % terminal

% dynamics on normalized time range [0,1]
% saturation input y normalized in [-1,1]
K = -[kp kd]/L;
y1 = K*x1; f1 = T*[x1(2); L*y1/I]; % linear regime
y2 = K*x2; f2 = T*[x2(2); L/I]; % upper sat
y3 = K*x3; f3 = T*[x3(2); -L/I]; % lower set

% test functions for each measure = monomials
g1 = mmon(x1,r); g2 = mmon(x2,r); g3 = mmon(x3,r);
g0 = mmon(x0,r); gT = mmon(xT,r);

% unknown moments of initial measure
y0 = mom(g0);

% unknown moments of terminal measure
yT = mom(gT);

% input LMI moment problem
cost = mom(xT'*xT);
Ay = mom(diff(g1,x1)*f1)+...
     mom(diff(g2,x2)*f2)+...
     mom(diff(g3,x3)*f3); % dynamics
% trajectory constraints
X = [y1^2<=1; y2>=1; y3<=-1]; 
% initial constraints
X0 = [x0(1)^2<=thetamax^2, x0(2)^2<=omegamax^2];
% terminal constraints
XT = [xT'*xT<=epsilon^2];
% bounds on trajectory
B = [x1'*x1<=4; x2'*x2<=4; x3'*x3<=4]; 

% input LMI moment problem
P = msdp(max(cost), ...
         mass(m1)+mass(m2)+mass(m3)==1, ...
         mass(m0)==1, ...
         Ay==yT-y0, ...
         X, X0, XT, B);

% solve LMI moment problem
[status,obj] = msol(P)

\end{verbatim}

With the help of this script and 
the SDP solver SeDuMi, we obtain the following sequence of upper bounds
(since we maximize)
on the maximum squared Euclidean norm of the final state:
\begin{center}
\begin{tabular}{c|cccc}
relaxation order $r$ & 1 & 2 & 3 & 4 \\\hline
upper bound & $1.0\cdot10^{-5}$ & $1.0\cdot10^{-5}$ & $1.0\cdot10^{-5}$ & $1.0\cdot10^{-5}$ \\
CPU time (sec.) & 0.2 & 0.5 & 0.7 & 0.9 \\
number of moments & 30 & 75 & 140 & 225
\end{tabular}
\end{center}
In the table we also indicate the CPU time (in seconds, on a standard desktop computer)
and the total number of moments (size of vector $y$ in the LMI relaxation).
We see that the bound obtained at the first relaxation ($r=1$)
is not modified for higher relaxations.
This clearly indicates that all initial conditions are
captured in the deadzone region at time $T$, which is
the box $[-2,2]\frac{10^{-1}\pi}{180} \times [-5,5]\frac{10^{-2}\pi}{180}
\supset \{x \in {\mathbb R}^2 \: :\: x^T x \leq 10^{-5}\}$.

\begin{figure}
\begin{center}
\includegraphics[width=0.5\textwidth]{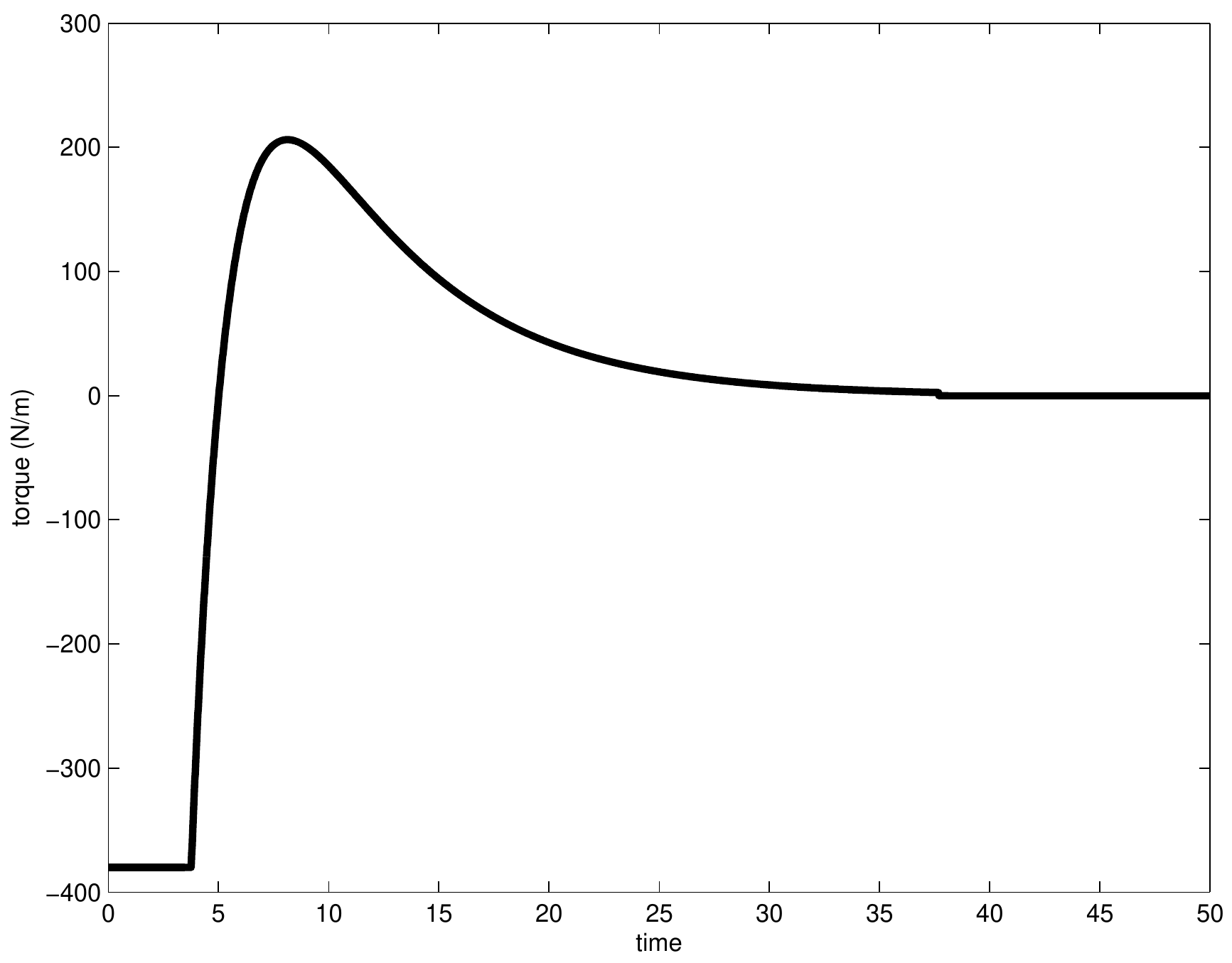}   
\caption{Torque input with lower saturation
during approx. 7\% of the time range.\label{fig:satur}}                              
\end{center}      
\end{figure}

If we want to use this approach to simulate a particular trajectory,
in the code we must modify the definition of the initial measure.
For example for initial conditions $x_1(0) = 50$, $x_2(0) = -1$,
we must insert the following sequence:
\begin{verbatim}
% given moments of initial measure = Dirac at x0
p = genpow(3,r); p = p(:,2:end); % powers
theta0 = 50; omega0 = -1; % in degrees
y0 = ones(size(p,1),1)*[theta0 omega0]*pi/180;
y0 = prod(y0.^p,2);
\end{verbatim}
As previously, the sequence of bounds on the maximum squared Euclidean norm 
of the final state is constantly equal to $1.0\cdot10^{-5}$, and in the
following table we represent as functions of the relaxation order $r$
the masses of measures $\mu_k$, $k=1,2,3$
which are indicators of the time spent by the trajectory in the
respective linear, upper saturation and lower saturation regimes:
\begin{center}
\begin{tabular}{c|cccccccc}
relaxation order $r$ & 1 & 2 & 3 & 4 & 5 & 6 & 7 \\ \hline
$\int d\mu_1$ & 37 & 89 & 92 & 92 & 93 & 93 & 93\\ 
$\int d\mu_2$ & 32 & 5.3 & 0.74 & 0.30 & 0.21 & 0.15 & 0.17 \\
$\int d\mu_3$ & 32 & 5.1 & 7.1 & 6.9 & 6.8 & 6.9 & 7.0
\end{tabular}
\end{center}
This indicates that most of the time (approx. 93\%) is spent in
the linear regime, with approx. 7\% of the time spent in the lower
saturation regime, and a negligible amount of time is spent in the
upper saturation regime. This is confirmed by simulation, see Figure \ref{fig:satur}.

\section{Notes and references}

Historically, the idea of reformulating nonconvex nonlinear ordinary differential equations (ODE) into 
convex LP, and especially linear partial differential equations (PDE)
in the space of probability measures, can be tracked back to the early 19th century.
It was Joseph Liouville in 1838 who first introduced the linear PDE involving the Jacobian
of the transformation exerted by the solution of an ODE on its initial condition \cite{l38}.
The idea was then largely expanded in Henri Poincar\'e's work on dynamical systems
at the end of the 19th century, see in particular 
\cite[Chapitre XII (Invariants int\'egraux)]{p99}. This work was pursued in the 20th
century in \cite{kb37}, \cite[Chapter VI (Systems with an integral invariant)]{ns47}
and more recently in the context of optimal transport by e.g. \cite{rr98}, \cite{v03}
or \cite{ags08}. The proof of Proposition \ref{continuity} can be found e.g.
in \cite[Chapter 8]{ags08}.
Poincar\'e himself in \cite[Section IV]{p08} mentions the potential of formulating
nonlinear ODEs as linear PDEs, and this programme has been carried out to
some extent by \cite{c32}, see also \cite{lm85}, \cite{ks91},
\cite{hl03}.
Our contribution is to apply the approach described in \cite{lhpt08}, see also
\cite{gq09}, to address
polynomial trajectory optimization problems. The use of LMI and measures
was also investigated in \cite{pr07} for building Lyapunov barrier
certificates, and based on a dual to Lyapunov's theorem described
in \cite{r01}. Our approach is similar, in the sense that optimization over systems
trajectories is formulated as an LP in the infinite-dimensional space
of measures. This LP problem is then approached as a generalized
moment problem via a hierarchy of LMI relaxations, following
the strategy described extensively in \cite{mom09}.
Finally, our control law validation problem of Sections \ref{valid}
and \ref{valid2} is comprehensively described in \cite{safev}.


\chapter{Polynomial optimal control}\label{chap:poc}

Our general setup for an {\bf optimal control problem} is the following:
\begin{equation}\label{ocp}
\begin{array}{lll}
p^* := & \inf & \int_0^T l(t,x(t),u(t))dt\\
& \mathrm{s.t.} &  \dot{x}(t) = f(t,x(t),u(t)),\\
& & x(t) \in X, \:\: u(t) \in U, \:\: t \in [0,T],\\
& & x(0) \in X_0, \quad x(T) \in X_T
\end{array}
\end{equation}
where the infimum is with respect to a control law $u : [0,T] \to {\mathbb R}^m$
which is a measurable function of time with values constrained to a given
set $U \subset {\mathbb R}^m$, and such that the resulting state trajectory
\[
x(t\:|\:x_0,u) = x_0 + \int_0^t f(s,x(s),u(s))ds
\]
starts at $x(0)=x_0$ in a given set $X_0$, terminates at time $T>0$ in
a given set $X_T$, and stays in a given set $X$ in between. It is
assumed that the given dynamics $f$ is smooth, so that there is a unique trajectory
given $x_0$ and $u$, which motivates our notation $x(t\:|\:x_0,u)$.
Also given is a smooth Lagrangian $l$.
The terminal time $T$ is either given or free.

\section{Controlled occupation measures}

As in Chapter \ref{chap:ido} we use occupation measures to model
problem (\ref{ocp}) as an infinite-dimensional LP. The
main difference however is that the occupation measures will
now depend on the control.

\begin{definition}[Controlled occupation measure]
Given an initial condition $x_0$ and a control law $u(t)$, the
controlled occupation measure of a trajectory $x(t\:|\:x_0,u)$
is defined as
\[
\mu(A\times B\times C\:|\:x_0,u) := \int_A I_B(x(t\:|\:x_0,u))dt
\]
for all $A \in {\mathscr B}([0,T])$, $B \in {\mathscr B}(X)$
and $C \in {\mathscr B}(U)$.
\end{definition}

A geometric interpretation is that $\mu$ measures the time spent
by the graph of the trajectory $(t,\:x(t\:|\:x_0,u),\:u(t))$
in a given subset $A\times B\times C$ of $[0,T]\times X\times U$.
An analytic interpretation is that integration w.r.t $\mu$ 
is equivalent to time-integration along a system trajectory.

If the initial condition $x_0 \in X$ is not a vector, but
an initial probability measure $\xi_0 \in {\mathscr P}(X)$, see
Definition \ref{initial}, we can proceed as in Section \ref{occmea}
and model the whole flow of trajectories with a measure.

\begin{definition}[Average controlled occupation measure]
Given an initial measure $\xi_0$ and a control law $u(t)$,
the average controlled occupation measure
of the flow of trajectories is defined as
\[
\mu(A\times B\times C\:|\:u) := \int_X \mu(A\times B\times C\:|\:x_0,u)\xi_0(dx_0)
\]
for all $A \in {\mathscr B}([0,T])$, $B \in {\mathscr B}(X)$
and $C \in {\mathscr B}(U)$.
\end{definition}

We also use the terminal measure $\xi_T$
as in Definition \ref{terminal}, and let 
$\mu_0 := \delta_0 \: \xi_0$, $\mu_T := \delta_T \: \xi_T$.
Measures $\mu$, $\mu_0$ and $\mu_T$ are linked by a linear PDE.
Let us now derive this equation
with the help of test functions $v$ depending on $t$ and $x$ only.
There is no dependence of $v$ on the control variable $u$
since the control law is an unknown in optimal
control problem (\ref{ocp}). 

Define the linear operator ${\mathcal L} : {\mathscr C}^1([0,T]\times X)
\to {\mathscr C}([0,T]\times X\times U)$ by
\[
v \mapsto {\mathcal L}v := \frac{\partial v}{\partial t} + (\mathrm{grad}\:v)' f
\]
and its adjoint operator ${\mathcal L}' : {\mathscr C}([0,T]\times X\times U)' \to
{\mathscr C}^1([0,T]\times X)'$ by
\[
\mu \mapsto {\mathcal L}'\mu =  -\frac{\partial \mu}{\partial t} - \mathrm{div}\:f\mu.
\]
Given a test function $v \in {\mathscr C}^1([0,T]\times X)$, it holds
\[
\begin{array}{rcl}
v(T,x(T)) & = & v(0,x(0)) + \int_0^T \dot{v}(t,x(t\:|\:x_0,u)dt \\
& = & v(0,x(0)) + \int_0^T {\mathcal L}v(t,x((t\:|\:x_0,u),u(t))dt \\
& = & v(0,x(0)) + \int_0^T \int_X \int_U {\mathcal L}v(t,x,u)\mu(dt,dx,du\:|\:x_0,u)
\end{array}
\]
and integrating w.r.t. $\xi_0$ we obtain
$\int {\mathcal L}v \mu = \int v\mu_T - \int v\mu_0$
for all $v$, which can be written in the sense of distributions as
${\mathcal L}'\mu = \mu_T - \mu_0$
or more explicitly
\begin{equation}\label{control}
\frac{\partial \mu}{\partial t} +
\mathrm{div}\:(f\mu) = \mu_0 - \mu_T.
\end{equation}
This is the {\bf controlled Liouville equation}.
The difference with the uncontrolled Liouville equation 
(\ref{liouville}) is that both $\mu$ and $f$ now also depend on the control variable $u$.
An occupation measure $\mu$ satisfying equation (\ref{control}) encodes
state trajectories but also control trajectories.

\section{Relaxed control}

We can disintegrate the occupation measure as
\begin{equation}\label{disint}
\mu(dt,dx,du) = dt\:\xi(dx\:|\:t)\:\omega(du\:|\:t,x)
\end{equation}
where the three components are as follows:
\begin{itemize}
\item $dt$ is the time marginal, the Lebesgue measure of time, corresponding to
the property that time flows uniformly;
\item $\xi(dx\:|\:t) \in {\mathscr P}(X)$ is the distribution of state conditional on $t$,
or state stochastic kernel,
a probability measure on $X$ for each $t \in [0,T]$, which models the
state interpreted as a time-dependent random variable;
\item $\omega(du\:|\:t,x) \in {\mathscr P}(U)$ is the distribution of the control conditional
on $t$ and $x$, or control stochastic kernel,
a probability measure on $U$ for each $t \in [0,T]$ and $x \in X$, which models
the control interpreted as a time- and state-dependent random variable.
\end{itemize}
It means that instead of a control law $u$ which is a measurable function of
time in $[0,T]$ with values in $U$, we have a {\bf relaxed control},
a probability measure
\[
\omega \in {\mathscr P}(U)
\]
parametrized in time $t \in [0,T]$ and space $x \in X$. Such parametrized
probability measures are called Young measures in the calculus of variations
and PDE literature. Our control, originally chosen as a measurable function
(of time and state), is therefore relaxed to a probability measure (parametrized in
time and state). Observe that the space of probability measures is larger
than any Lebesgue space, since for the particular choice of a time-dependent
Dirac measure
\[
\omega(du\:|\:t,x) = \delta_{u(t,x)}
\]
with $u(t,x) \in U$
we retrieve a classical control law which is a function of time and state.

The controlled Liouville equation (\ref{control}) can be written as
\[
\begin{array}{rcl}
\int v \mu_T - \int v \mu_0 & = & \int {\mathcal L}v \mu \\
& = & \int_T\int_X\int_U \left(\frac{{\partial v}(t,x)}{\partial t}
+ (\mathrm{grad}\:v(t,x))'f(t,x,u)\right)\omega(du\:|\:t,x)\xi(dx\:|\:t)dt \\
& = & \int_T\int_X \left(\frac{{\partial v}(t,x)}{\partial t}
+ (\mathrm{grad}\:v(t,x))'\left[\int_U f(t,x,u)\omega(du\:|\:t,x)\right]
\right)\xi(dx\:|\:t)dt
\end{array}
\]
for all test functions $v \in {\mathscr C}^1([0,T]\times X)$.
It is now apparent that the trajectories modeled by the
controlled Liouville equation are generated by a family of absolutely
continuous trajectories of the relaxed controlled ODE
\[
\dot{x}(t) = \int_U f(t,x(t),u)\omega(du\:|\:t,x).
\]
Indeed, in optimal control problem (\ref{ocp}), the original control system
\begin{equation}\label{code}
\dot{x}(t) = f(t,x(t),u(t)), \quad u(t) \in U
\end{equation}
can be interpreted as a differential inclusion
\[
\dot{x}(t) \in f(t,x(t),U) := \{f(t,x(t),u) \: :\: u \in U\}
\]
where the state velocity vector $\dot{x}(t)$ can be chosen
anywhere in the set $f(t,x(t),U) \subset {\mathbb R}^n$.
In contrast, any triplet of measures $(\mu,\mu_0,\mu_T)$ satisfying
the controlled Liouville equation (\ref{control}) corresponds
to a family of trajectories of the relaxed, or convexified differential inclusion
\[
\dot{x}(t) \in \mathrm{conv}\:f(t,x(t),U).
\]
In that sense, the set of trajectories modeled by the
controlled Liouville equation (\ref{control}) is
larger than the set of trajectories of the control
system (\ref{code}). As will be seen in the numerical
example section, this is an advantage of the occupation
measure framework, in the sense that we will be able
to construct relaxed or stochastic control laws
that cannot can obtained using functions.

Based on the above discussion, we can define the
{\bf relaxed optimal control problem}
\begin{equation}\label{rocp}
\begin{array}{lll}
p^*_R := & \inf & \int_0^T l(t,x(t),u(t))dt\\
& \mathrm{s.t.} &  \dot{x}(t) \in \mathrm{conv}\:f(t,x(t),U),\\
& & x(t) \in X, \:\: u(t) \in U, \:\: t \in [0,T],\\
& & x(0) \in X_0, \quad x(T) \in X_T
\end{array}
\end{equation}
and it holds $p^*_R\leq p^*$. Contrived optimal control problems
(e.g. with stringent state constraints)
can be cooked up such that $p^*_R<p^*$, but generically (in a
sense to be defined rigorously, but not in this document),
the following assumption will be satisfied.

\begin{assumption}[No relaxation gap]\label{nogap}
We assume that $p^*_R=p^*$.
\end{assumption}

\section{Measure LP}

Using the controlled occupation measure and relaxed controls
of the previous sections, and under Assumption \ref{nogap},
relaxed optimal control problem 
(\ref{rocp}) can be
formulated as an infinite-dimensional measure LP
\[
\begin{array}{lll}
p^* = & \inf & \langle l,\mu \rangle \\
& \mathrm{s.t.} & \frac{\partial \mu}{\partial t} + \mathrm{div}\:f\mu = \mu_0 - \mu_T
\end{array}
\]
where the infimum is w.r.t. measures $(\mu,\:\mu_0,\:\mu_T) \in
{\mathscr M}_+([0,T]\times X\times U)\times{\mathscr M}_+(\{0\}\times X_0))
\times{\mathscr M}_+(\{T\}\times X_T)$.
We can then rely on the results of Section \ref{inflplmi}
to build a hierarchy of finite-dimensional
LMI relaxations for this problem. This generates a monotonically
nondecreasing sequence of lower bounds asymptotically converging
to $p^*$. Details are omitted.

\begin{remark}[Autonomous case]\label{autonomousocp}
If the terminal time $T$ is free and the Lagrangian $l$ and the dynamics
$f$ do not depend explicitly on time $t$, then the measure LP becomes
\[
\begin{array}{lll}
p^* = & \inf & \langle l,\mu \rangle \\
& \mathrm{s.t.} & \mathrm{div}\:f\mu = \mu_0 - \mu_T\\
\end{array}
\]
where the minimum is now taken w.r.t.
$(\mu,\:\mu_0,\:\mu_T) \in {\mathscr M}_+(X\times U)\times
{\mathscr M}_+(X_0)\times{\mathscr M}_+(X_T)$, i.e.
the measures do not depend explicitly on time either.
\end{remark}

\begin{example}[Linear quadratic regulator]\label{ex:lqr}
Consider the elementary scalar linear quadratic regulator problem
\[
\begin{array}{lll}
p^* = & \inf & \int_0^T (x^2(t)+u^2(t))dt\\
& \mathrm{s.t.} &  \dot{x}(t) = u(t), \quad t\in[0,T]\\	
& & x(0) = 1, \quad x(T) = 0
\end{array}
\]
with given initial and terminal conditions.
The corresponding autonomous measure LP is
\[
\begin{array}{lll}
p^* = & \inf & \langle x^2+u^2,\mu \rangle\\
& \mathrm{s.t.} &  \mathrm{div}\:u\:\mu = \delta_{1} - \delta_{0}\\
\end{array}
\]
where the minimum is w.r.t. occupation measure $\mu$.
Its moment LP problem reads
\[
\begin{array}{lll}
p^* = & \inf & \int (x^2+u^2)\mu(dx,du) \\
& \mathrm{s.t.} &  \alpha \int x^{\alpha-1}u\mu(dx,du) =  -1, \quad \alpha=0,1,2\ldots
\end{array}
\]
or equivalently
\[
\begin{array}{lll}
p^* = & \inf & y_{20}+y_{02} \\
& \mathrm{s.t.} & y_{01} = 2y_{11} = 3y_{21} = \cdots = -1 \\
&& M(y) \geq 0
\end{array}
\]
where the minimum is w.r.t. the moments of the occupation measure 
\[
y_{\alpha} = \int x^{\alpha_1}u^{\alpha_2}\mu(dx,du), \quad \alpha=0,1,2\ldots
\]
Solving the first LMI relaxation
\[
\begin{array}{lll}
p^*_1 = & \inf & y_{20}+y_{02} \\
& \mathrm{s.t.} & y_{01} = 2y_{11} = -1 \\
&& M_1(y) = \left(\begin{array}{ccc}
y_{00} & y_{10} & y_{01} \\
y_{10} & y_{20} & y_{11} \\
y_{01} & y_{11} & y_{02} 
\end{array}\right) \geq 0.
\end{array}
\]
with the SDP solver SeDuMi yields (rounded to 3 significant digits)
\[
M_1(y^*) = \left(\begin{array}{rrr}
3.66 & 1.00 & -1.00 \\
1.00 & 0.500 & -0.500\\
-1.00  & -0.500 & 0.500
\end{array}\right).
\]
This example can be solved analytically (with a scalar Riccati equation)
and the solution is the state-feedback $u(t)=-x(t)$
corresponding to the optimal trajectory $x(t)=e^{-t}$
with cost $p^*=\int_0^{\infty} 2e^{-2t}dt = 1$ and the
optimal occupation measure
\[
\mu(dx,du) = \int_0^{\infty} \delta_{e^{-t}}(dx) \delta_{-e^{-t}}(du)dt
\]
with moments
\[
y_{\alpha} = (-1)^{\alpha_2} \int_0^{\infty} e^{-(\alpha_1+\alpha_2)t}dt, \quad
\alpha \in {\mathbb N}^2
\]
equal to $y_{00} = \infty$, $y_{10} = 1$, $y_{01} = -1$,
$y_{20} = \frac{1}{2}$, $y_{11} = -\frac{1}{2}$, $y_{02} = \frac{1}{2}$ etc.
We observe that the numerical moments $y^*$ closely match, except the
mass $y_{00}$ which should approximate the terminal time $T$.
Note however that for this numerical value of $T\approx 3.66$, the cost is
$\int_0^T 2e^{-2t}dt \approx 0.999$ almost equal to the optimal value
$p^*=1$.
\end{example}

\section{Optimal control recovery}

Once an LMI relaxation of given order is solved,
we expect vector $y$ to contain approximate moments
of the optimal occupation measure corresponding to
the optimal trajectory (if it is unique),
or at least a superposition (convex combination)
of optimal trajectories. To recover the optimal control,
or the optimal state trajectory, we can use the dual
problem, which is a relaxation of the Hamilton-Jacobi-Bellman
PDE of optimal control. However, we do not elaborate further
on these techniques in this document.

\section{Back to the motivating example}\label{bolza2}

Let us come back to Bolza's example of Section \ref{bolza}.
We saw that the infimum can be approached by a control sequence switching
increasingly quickly between $-1$ and $+1$, so the idea is to relax the ODE
\[
\dot{x}(t) = f(t,x(t),u(t))
\]
with the following differential equation
\[
\dot{x}(t) = \int f(t,x(t),u)\omega(du\:|\:t)
\]
where $\omega(du | t)$ is a probability measure
parametrized in $t$. 
State trajectories are then obtained by integration w.r.t. time
and control
\[
x(t) = x(0) + \int_0^t \int f(s,x(s),u)\omega(du\:|\:s)ds.
\]
Here for the Bolza example we choose
\[
\omega(du\:|\:t) = \frac{1}{2}(\delta_{u=-1} + \delta_{u=+1})
\]
a time-independent weighted sum of two Dirac measures
at $u=-1$ and $u=+1$. The relaxed state trajectory is then
equal to
\[
\begin{array}{rcl}
x(t) & = & \frac{1}{2}\left(\int_0^t f(s,x(s),-1)ds +
\int_0^t f(s,x(s),+1)ds\right) \\
& = & \frac{1}{2}\left(-\int_0^t ds + \int_0^t ds\right) = 0
\end{array}
\]
and the relaxed objective function is equal to
\[
\begin{array}{rcl} 
\int_0^1 \int_U l(t,x(t),u)d\omega(u\:|\:t)dt & = & \frac{1}{2}\left(\int_0^1 l(t,x(t),-1)
+ \int_0^1 l(t,x(t),+1)\right) \\
& = & \int_0^1 x^4(t)dt = 0
\end{array}
\]
so that the infimum $p^*=0$ is reached.

The corresponding GloptiPoly script is as follows:

\begin{verbatim}
% initial point measure
mpol('t0'); mpol('x0');

% occupation measure
mpol('t'); mpol('x'); mpol('u');
meas(t,2); meas(x,2); meas(u,2);

% final point measure
mpol('tT'); mpol('xT',1);
meas(tT,3); meas(xT,3);

r = input('order of relaxation ='); r = 2*r;

% define test function arrays
v = mmon([t;x],r);
v0 = mmon([t0; x0],r);
vT = mmon([tT; xT],r);

% dynamics f(x,u) = u
dvdt = diff(v,t) + diff(v,x)*u;

% assign initial point
assign([t0; x0],[ 0; 0]);

% input LMI moment problem
P = msdp (min(x^2 + (1-u^2)^2 ), ...
   0 == mom(dvdt) + double(v0) - mom(vT), ...
   t*(1-t) >= 0, u^2 <= 1, x^2 <= 1, ...
   tT == 1, xT^2 <= 1);

% solve LMI moment problem
[status, obj]  =  msol(P)
\end{verbatim}

\section{Notes and references}

The use of relaxations and LP formulations of optimal control problems (on
ordinary differential equations and partial differential equations) is classical, and can be
traced back to the work by L. C. Young \cite{y69},
Filippov \cite{f67}, Warga \cite{w72}, Gamkrelidze \cite{g75}
and Rubio \cite{r86} amongst many others. For more details and a historical survey,
see e.g. \cite[Part III]{f99}.
Parametrized measures arising in the disintegration (\ref{disint}) of the
occupation measures are called Young measures in the PDE
literature, see e.g. \cite{p97} or \cite{r97}.
Our contribution is to notice that hierarchies of LMI relaxations
can be used to solve numerically the infinite-dimensional LP on measures
arising from relaxed optimal control problems with polynomial data,
following the methodology described originally in \cite{lhpt08}.
By the Filippov-Wa$\dot{\mathrm z}$ewski relaxation theorem \cite{af90},
the trajectories of the optimal control problem (\ref{ocp})
are dense (w.r.t. the metric of uniform convergence of absolutely continuous
functions of time) in the set of trajectories of the relaxed optimal control problem
(\ref{rocp}). This justifies Assumption \ref{nogap}.
Finally, our motivating Bolza problem of Sections \ref{bolza} and \ref{bolza2}
is a classical example of calculus of variations
illustrating that an optimal control problem with smooth data (infinitely differentiable
Lagrangian and dynamics, no state or input constraints) can have
a highly oscillatory optimal solution,
see e.g. \cite[Example 4.8]{d07}.

\end{document}